\documentclass[11pt]{amsart}
\usepackage{graphicx}
\usepackage{amssymb}
\usepackage{epsfig}
\usepackage{anysize}

\newcommand{\tlowername}[2]%
{$\stackrel{\makebox[1pt]{#1}}%
{\begin{picture}(0,0)%
\put(0,0){\makebox(0,6)[t]{\makebox[1pt]{$#2$}}}%
\end{picture}}$}%

\newcommand{\AR}[1]%
{\begin{picture}(#1,0)%
\put(0,0){\vector(1,0){#1}}%
\end{picture}}%

\newcommand{\DOTAR}[1]%
{\NUMBEROFDOTS=#1%
\divide\NUMBEROFDOTS by 3%
\begin{picture}(#1,0)%
\multiput(0,0)(3,0){\NUMBEROFDOTS}{\circle*{1}}%
\put(#1,0){\vector(1,0){0}}%
\end{picture}}%

\newcommand{\MONO}[1]%
{\begin{picture}(#1,0)%
\put(0,0){\vector(1,0){#1}}%
\put(2,-2){\line(0,1){4}}%
\end{picture}}%

\newcommand{\EPI}[1]%
{\begin{picture}(#1,0)(-#1,0)%
\put(-#1,0){\vector(1,0){#1}}%
\put(-6,-2){\line(0,1){4}}%
\end{picture}}%

\newcommand{\BIMO}[1]%
{\begin{picture}(#1,0)(-#1,0)%
\put(-#1,0){\vector(1,0){#1}}%
\put(-6,-2){\line(0,1){4}}%
\put(-#1,-2){\hspace{2pt}\line(0,1){4}}%
\end{picture}}%

\newcommand{\BIAR}[1]%
{\begin{picture}(#1,4)%
\put(0,0){\vector(1,0){#1}}%
\put(0,4){\vector(1,0){#1}}%
\end{picture}}%

\newcommand{\EQL}[1]%
{\begin{picture}(#1,0)%
\put(0,1){\line(1,0){#1}}%
\put(0,-1){\line(1,0){#1}}%
\end{picture}}%

\newcommand{\ADJAR}[1]%
{\begin{picture}(#1,4)%
\put(0,0){\vector(1,0){#1}}%
\put(#1,4){\vector(-1,0){#1}}%
\end{picture}}%

\newcommand{\BKAR}[1]%
{\begin{picture}(#1,0)%
\put(#1,0){\vector(-1,0){#1}}%
\end{picture}}%

\newcommand{\BKDOTAR}[1]%
{\NUMBEROFDOTS=#1%
\divide\NUMBEROFDOTS by 3%
\begin{picture}(#1,0)%
\multiput(#1,0)(-3,0){\NUMBEROFDOTS}{\circle*{1}}%
\put(0,0){\vector(-1,0){0}}%
\end{picture}}%

\newcommand{\BKMONO}[1]%
{\begin{picture}(#1,0)(-#1,0)%
\put(0,0){\vector(-1,0){#1}}%
\put(-2,-2){\line(0,1){4}}%
\end{picture}}%

\newcommand{\BKEPI}[1]%
{\begin{picture}(#1,0)%
\put(#1,0){\vector(-1,0){#1}}%
\put(6,-2){\line(0,1){4}}%
\end{picture}}%

\newcommand{\BKBIMO}[1]%
{\begin{picture}(#1,0)%
\put(#1,0){\vector(-1,0){#1}}%
\put(6,-2){\line(0,1){4}}%
\put(#1,-2){\hspace{-2pt}\line(0,1){4}}%
\end{picture}}%

\newcommand{\BKBIAR}[1]%
{\begin{picture}(#1,4)%
\put(#1,0){\vector(-1,0){#1}}%
\put(#1,4){\vector(-1,0){#1}}%
\end{picture}}%

\newcommand{\BKADJAR}[1]%
{\begin{picture}(#1,4)%
\put(0,4){\vector(1,0){#1}}%
\put(#1,0){\vector(-1,0){#1}}%
\end{picture}}%

\newcommand{\lowername}[2]%
{$\stackrel{\makebox[1pt]{#1}}%
{\begin{picture}(0,0)%
\truex{600}%
\put(0,0){\makebox(0,\value{x})[t]{\makebox[1pt]{$#2$}}}%
\end{picture}}$}%

\newcommand{\hcase}[2]%
{\makebox[0pt]%
{\raisebox{-1pt}[0pt][0pt]{#1{#2}}}}%

\newcommand{\Hcase}[3]%
{\makebox[0pt]
{\raisebox{-1pt}[0pt][0pt]%
{$\stackrel{\makebox[0pt]{$\textstyle{#2}$}}{#1{#3}}$}}}%

\newcommand{\hcasE}[3]%
{\makebox[0pt]%
{\raisebox{-9pt}[0pt][0pt]%
{\lowername{#1{#3}}{#2}}}}%

\newcommand{\hbicase}[2]%
{\makebox[0pt]%
{\raisebox{-2.5pt}[0pt][0pt]{#1{#2}}}}%

\newcommand{\Hbicase}[4]%
{\makebox[0pt]
{\raisebox{-10.5pt}[0pt][0pt]%
{$\stackrel{\makebox[0pt]{$\textstyle{#2}$}}%
{\mbox{\lowername{#1{#4}}{#3}}}$}}}%

\newcommand{\EAR}[1]%
{\begin{picture}(#1,0)%
\put(0,0){\vector(1,0){#1}}%
\end{picture}}%

\newcommand{\EDOTAR}[1]%
{\truex{100}\truey{300}%
\NUMBEROFDOTS=#1%
\divide\NUMBEROFDOTS by \value{y}%
\begin{picture}(#1,0)%
\multiput(0,0)(\value{y},0){\NUMBEROFDOTS}%
{\circle*{\value{x}}}%
\put(#1,0){\vector(1,0){0}}%
\end{picture}}%

\newcommand{\EMONO}[1]%
{\begin{picture}(#1,0)%
\put(0,0){\vector(1,0){#1}}%
\truex{300}\truey{600}%
\put(\value{x},-\value{x}){\line(0,1){\value{y}}}%
\end{picture}}%

\newcommand{\EEPI}[1]%
{\begin{picture}(#1,0)(-#1,0)%
\put(-#1,0){\vector(1,0){#1}}%
\truex{300}\truey{600}\truez{800}%
\put(-\value{z},-\value{x}){\line(0,1){\value{y}}}%
\end{picture}}%

\newcommand{\EBIMO}[1]%
{\begin{picture}(#1,0)(-#1,0)%
\put(-#1,0){\vector(1,0){#1}}%
\truex{300}\truey{600}\truez{800}%
\put(-\value{z},-\value{x}){\line(0,1){\value{y}}}%
\put(-#1,-\value{x}){\hspace{3pt}\line(0,1){\value{y}}}%
\end{picture}}%

\newcommand{\EBIAR}[1]%
{\truex{400}%
\begin{picture}(#1,\value{x})%
\put(0,0){\vector(1,0){#1}}%
\put(0,\value{x}){\vector(1,0){#1}}%
\end{picture}}%

\newcommand{\EEQL}[1]%
{\begin{picture}(#1,0)%
\truex{200}%
\put(0,\value{x}){\line(1,0){#1}}%
\put(0,0){\line(1,0){#1}}%
\end{picture}}%

\newcommand{\EADJAR}[1]%
{\truex{400}%
\begin{picture}(#1,\value{x})%
\put(0,0){\vector(1,0){#1}}%
\put(#1,\value{x}){\vector(-1,0){#1}}%
\end{picture}}%

\newcommand{\ear}%
{\hspace{\SOURCE\unitlength}%
\hcase{\EAR}{\ARROWLENGTH}}%

\newcommand{\Ear}[1]%
{\hspace{\SOURCE\unitlength}%
\Hcase{\EAR}{#1}{\ARROWLENGTH}}%

\newcommand{\eaR}[1]%
{\hspace{\SOURCE\unitlength}%
\hcasE{\EAR}{#1}{\ARROWLENGTH}}%

\newcommand{\edotar}%
{\hspace{\SOURCE\unitlength}%
\hcase{\EDOTAR}{\ARROWLENGTH}}%

\newcommand{\Edotar}[1]%
{\hspace{\SOURCE\unitlength}%
\Hcase{\EDOTAR}{#1}{\ARROWLENGTH}}%

\newcommand{\edotaR}[1]%
{\hspace{\SOURCE\unitlength}%
\hcasE{\EDOTAR}{#1}{\ARROWLENGTH}}%

\newcommand{\emono}%
{\hspace{\SOURCE\unitlength}%
\hcase{\EMONO}{\ARROWLENGTH}}%

\newcommand{\Emono}[1]%
{\hspace{\SOURCE\unitlength}%
\Hcase{\EMONO}{#1}{\ARROWLENGTH}}%

\newcommand{\emonO}[1]%
{\hspace{\SOURCE\unitlength}%
\hcasE{\EMONO}{#1}{\ARROWLENGTH}}%

\newcommand{\eepi}%
{\hspace{\SOURCE\unitlength}%
\hcase{\EEPI}{\ARROWLENGTH}}%

\newcommand{\Eepi}[1]%
{\hspace{\SOURCE\unitlength}%
\Hcase{\EEPI}{#1}{\ARROWLENGTH}}%

\newcommand{\eepI}[1]%
{\hspace{\SOURCE\unitlength}%
\hcasE{\EEPI}{#1}{\ARROWLENGTH}}%

\newcommand{\ebimo}%
{\hspace{\SOURCE\unitlength}%
\hcase{\EBIMO}{\ARROWLENGTH}}%

\newcommand{\Ebimo}[1]%
{\hspace{\SOURCE\unitlength}%
\Hcase{\EBIMO}{#1}{\ARROWLENGTH}}%

\newcommand{\ebimO}[1]%
{\hspace{\SOURCE\unitlength}%
\hcasE{\EBIMO}{#1}{\ARROWLENGTH}}%

\newcommand{\eiso}%
{\hspace{\SOURCE\unitlength}%
\Hcase{\EAR}{\cong}{\ARROWLENGTH}}%

\newcommand{\Eiso}[1]%
{\hspace{\SOURCE\unitlength}%
\Hcase{\EAR}{\cong#1}{\ARROWLENGTH}}%

\newcommand{\eisO}[1]%
{\hspace{\SOURCE\unitlength}%
\hcasE{\EAR}{\cong#1}{\ARROWLENGTH}}%

\newcommand{\ebiar}%
{\hspace{\SOURCE\unitlength}%
\hbicase{\EBIAR}{\ARROWLENGTH}}%

\newcommand{\Ebiar}[2]%
{\hspace{\SOURCE\unitlength}%
\Hbicase{\EBIAR}{#1}{#2}{\ARROWLENGTH}}%

\newcommand{\eeql}%
{\hspace{\SOURCE\unitlength}%
\hbicase{\EEQL}{\ARROWLENGTH}}%

\newcommand{\eadjar}%
{\hspace{\SOURCE\unitlength}%
\hbicase{\EADJAR}{\ARROWLENGTH}}%

\newcommand{\Eadjar}[2]%
{\hspace{\SOURCE\unitlength}%
\Hbicase{\EADJAR}{#1}{#2}{\ARROWLENGTH}}%

\newcommand{\WAR}[1]%
{\begin{picture}(#1,0)%
\put(#1,0){\vector(-1,0){#1}}%
\end{picture}}%

\newcommand{\WDOTAR}[1]%
{\truex{100}\truey{300}%
\NUMBEROFDOTS=#1%
\divide\NUMBEROFDOTS by \value{y}%
\begin{picture}(#1,0)%
\multiput(#1,0)(-\value{y},0){\NUMBEROFDOTS}%
{\circle*{\value{x}}}%
\put(0,0){\vector(-1,0){0}}%
\end{picture}}%

\newcommand{\WMONO}[1]%
{\begin{picture}(#1,0)(-#1,0)%
\put(0,0){\vector(-1,0){#1}}%
\truex{300}\truey{600}%
\put(-\value{x},-\value{x}){\line(0,1){\value{y}}}%
\end{picture}}%

\newcommand{\WEPI}[1]%
{\begin{picture}(#1,0)%
\put(#1,0){\vector(-1,0){#1}}%
\truex{300}\truey{600}\truez{800}%
\put(\value{z},-\value{x}){\line(0,1){\value{y}}}%
\end{picture}}%

\newcommand{\WBIMO}[1]%
{\begin{picture}(#1,0)%
\put(#1,0){\vector(-1,0){#1}}%
\truex{300}\truey{600}\truez{800}%
\put(\value{z},-\value{x}){\line(0,1){\value{y}}}%
\put(#1,-\value{x}){\hspace{-3pt}\line(0,1){\value{y}}}%
\end{picture}}%

\newcommand{\WBIAR}[1]%
{\truex{400}%
\begin{picture}(#1,\value{x})%
\put(#1,0){\vector(-1,0){#1}}%
\put(#1,\value{x}){\vector(-1,0){#1}}%
\end{picture}}%

\newcommand{\WADJAR}[1]%
{\truex{400}%
\begin{picture}(#1,\value{x})%
\put(0,\value{x}){\vector(1,0){#1}}%
\put(#1,0){\vector(-1,0){#1}}%
\end{picture}}%

\newcommand{\war}%
{\hspace{\SOURCE\unitlength}%
\hcase{\WAR}{\ARROWLENGTH}}%

\newcommand{\War}[1]%
{\hspace{\SOURCE\unitlength}%
\Hcase{\WAR}{#1}{\ARROWLENGTH}}%

\newcommand{\waR}[1]%
{\hspace{\SOURCE\unitlength}%
\hcasE{\WAR}{#1}{\ARROWLENGTH}}%

\newcommand{\wdotar}%
{\hspace{\SOURCE\unitlength}%
\hcase{\WDOTAR}{\ARROWLENGTH}}%

\newcommand{\Wdotar}[1]%
{\hspace{\SOURCE\unitlength}%
\Hcase{\WDOTAR}{#1}{\ARROWLENGTH}}%

\newcommand{\wdotaR}[1]%
{\hspace{\SOURCE\unitlength}%
\hcasE{\WDOTAR}{#1}{\ARROWLENGTH}}%

\newcommand{\wmono}%
{\hspace{\SOURCE\unitlength}%
\hcase{\WMONO}{\ARROWLENGTH}}%

\newcommand{\Wmono}[1]%
{\hspace{\SOURCE\unitlength}%
\Hcase{\WMONO}{#1}{\ARROWLENGTH}}%

\newcommand{\wmonO}[1]%
{\hspace{\SOURCE\unitlength}%
\hcasE{\WMONO}{#1}{\ARROWLENGTH}}%

\newcommand{\wepi}%
{\hspace{\SOURCE\unitlength}%
\hcase{\WEPI}{\ARROWLENGTH}}%

\newcommand{\Wepi}[1]%
{\hspace{\SOURCE\unitlength}%
\Hcase{\WEPI}{#1}{\ARROWLENGTH}}%

\newcommand{\wepI}[1]%
{\hspace{\SOURCE\unitlength}%
\hcasE{\WEPI}{#1}{\ARROWLENGTH}}%

\newcommand{\wbimo}%
{\hspace{\SOURCE\unitlength}%
\hcase{\WBIMO}{\ARROWLENGTH}}%

\newcommand{\Wbimo}[1]%
{\hspace{\SOURCE\unitlength}%
\Hcase{\WBIMO}{#1}{\ARROWLENGTH}}%

\newcommand{\wbimO}[1]%
{\hspace{\SOURCE\unitlength}%
\hcasE{\WBIMO}{#1}{\ARROWLENGTH}}%

\newcommand{\wiso}%
{\hspace{\SOURCE\unitlength}%
\Hcase{\WAR}{\cong}{\ARROWLENGTH}}%

\newcommand{\Wiso}[1]%
{\hspace{\SOURCE\unitlength}%
\Hcase{\WAR}{#1}{\ARROWLENGTH}}%

\newcommand{\wisO}[1]%
{\hspace{\SOURCE\unitlength}%
\hcasE{\WAR}{#1}{\ARROWLENGTH}}%

\newcommand{\wbiar}%
{\hspace{\SOURCE\unitlength}%
\hbicase{\WBIAR}{\ARROWLENGTH}}%

\newcommand{\Wbiar}[2]%
{\hspace{\SOURCE\unitlength}%
\Hbicase{\WBIAR}{#1}{#2}{\ARROWLENGTH}}%

\newcommand{\weql}%
{\hspace{\SOURCE\unitlength}%
\hbicase{\EEQL}{\ARROWLENGTH}}%

\newcommand{\wadjar}%
{\hspace{\SOURCE\unitlength}%
\hbicase{\WADJAR}{\ARROWLENGTH}}%

\newcommand{\Wadjar}[2]%
{\hspace{\SOURCE\unitlength}%
\Hbicase{\WADJAR}{#1}{#2}{\ARROWLENGTH}}%

\newcommand{\vcase}[2]{#1{#2}}%

\newcommand{\Vcase}[3]{\makebox[0pt]%
{\makebox[0pt][r]{\raisebox{0pt}[0pt][0pt]{${#2}\hspace{2pt}$}}}#1{#3}}%

\newcommand{\vcasE}[3]{\makebox[0pt]%
{#1{#3}\makebox[0pt][l]{\raisebox{0pt}[0pt][0pt]{\hspace{2pt}$#2$}}}}%

\newcommand{\vbicase}[2]{\makebox[0pt]{{#1{#2}}}}%

\newcommand{\Vbicase}[4]{\makebox[0pt]%
{\makebox[0pt][r]{\raisebox{0pt}[0pt][0pt]{$#2$\hspace{4pt}}}#1{#4}%
\makebox[0pt][l]{\raisebox{0pt}[0pt][0pt]{\hspace{5pt}$#3$}}}}%

\newcommand{\SAR}[1]%
{\begin{picture}(0,0)%
\put(0,0){\makebox(0,0)%
{\begin{picture}(0,#1)%
\put(0,#1){\vector(0,-1){#1}}%
\end{picture}}}\end{picture}}%

\newcommand{\SDOTAR}[1]%
{\truex{100}\truey{300}%
\NUMBEROFDOTS=#1%
\divide\NUMBEROFDOTS by \value{y}%
\begin{picture}(0,0)%
\put(0,0){\makebox(0,0)%
{\begin{picture}(0,#1)%
\multiput(0,#1)(0,-\value{y}){\NUMBEROFDOTS}%
{\circle*{\value{x}}}%
\put(0,0){\vector(0,-1){0}}%
\end{picture}}}\end{picture}}%

\newcommand{\SMONO}[1]%
{\begin{picture}(0,0)%
\put(0,0){\makebox(0,0)%
{\begin{picture}(0,#1)%
\put(0,#1){\vector(0,-1){#1}}%
\truex{300}\truey{600}%
\put(0,#1){\begin{picture}(0,0)%
\put(-\value{x},-\value{x}){\line(1,0){\value{y}}}\end{picture}}%
\end{picture}}}\end{picture}}%

\newcommand{\SEPI}[1]%
{\begin{picture}(0,0)%
\put(0,0){\makebox(0,0)%
{\begin{picture}(0,#1)%
\put(0,#1){\vector(0,-1){#1}}%
\truex{300}\truey{600}\truez{800}%
\put(-\value{x},\value{z}){\line(1,0){\value{y}}}%
\end{picture}}}\end{picture}}%

\newcommand{\SBIMO}[1]%
{\begin{picture}(0,0)%
\put(0,0){\makebox(0,0)%
{\begin{picture}(0,#1)%
\put(0,#1){\vector(0,-1){#1}}%
\truex{300}\truey{600}\truez{800}%
\put(0,#1){\begin{picture}(0,0)%
\put(-\value{x},-\value{x}){\line(1,0){\value{y}}}\end{picture}}%
\put(-\value{x},\value{z}){\line(1,0){\value{y}}}%
\end{picture}}}\end{picture}}%

\newcommand{\SBIAR}[1]%
{\begin{picture}(0,0)%
\truex{200}%
\put(0,0){\makebox(0,0)%
{\begin{picture}(0,#1)\put(-\value{x},#1){\vector(0,-1){#1}}%
\put(\value{x},#1){\vector(0,-1){#1}}%
\end{picture}}}\end{picture}}%

\newcommand{\SEQL}[1]%
{\begin{picture}(0,0)%
\truex{100}%
\put(0,0){\makebox(0,0)%
{\begin{picture}(0,#1)\put(-\value{x},#1){\line(0,-1){#1}}%
\put(\value{x},#1){\line(0,-1){#1}}%
\end{picture}}}\end{picture}}%

\newcommand{\sarv}[1]{\vcase{\SAR}{#100}}%

\newcommand{\sar}{\sarv{50}}%

\newcommand{\Sarv}[2]{\Vcase{\SAR}{#1}{#200}}%

\newcommand{\Sar}[1]{\Sarv{#1}{50}}%

\newcommand{\saRv}[2]{\vcasE{\SAR}{#1}{#200}}%

\newcommand{\saR}[1]{\saRv{#1}{50}}%

\newcommand{\Sisov}[2]%
{\Vbicase{\SAR}{#1\hspace{-2pt}}{\hspace{-2pt}\cong}{#200}}%

\newcommand{\seqlv}[1]{\vbicase{\SEQL}{#100}}%

\newcommand{\seql}{\seqlv{50}}%

\newcommand{\NAR}[1]%
{\begin{picture}(0,0)%
\put(0,0){\makebox(0,0)%
{\begin{picture}(0,#1)\put(0,0){\vector(0,1){#1}}%
\end{picture}}}\end{picture}}%

\newcommand{\NDOTAR}[1]%
{\truex{100}\truey{300}%
\NUMBEROFDOTS=#1%
\divide\NUMBEROFDOTS by \value{y}%
\begin{picture}(0,0)%
\put(0,0){\makebox(0,0)%
{\begin{picture}(0,#1)%
\multiput(0,0)(0,\value{y}){\NUMBEROFDOTS}%
{\circle*{\value{x}}}%
\put(0,#1){\vector(0,1){0}}%
\end{picture}}}\end{picture}}%

\newcommand{\NMONO}[1]%
{\begin{picture}(0,0)%
\put(0,0){\makebox(0,0)%
{\begin{picture}(0,#1)%
\put(0,0){\vector(0,1){#1}}%
\truex{300}\truey{600}%
\put(-\value{x},\value{x}){\line(1,0){\value{y}}}%
\end{picture}}}%
\end{picture}}%

\newcommand{\NEPI}[1]%
{\begin{picture}(0,0)%
\put(0,0){\makebox(0,0)%
{\begin{picture}(0,#1)%
\put(0,0){\vector(0,1){#1}}%
\truex{300}\truey{600}\truez{800}%
\put(0,#1){\begin{picture}(0,0)%
\put(-\value{x},-\value{z}){\line(1,0){\value{y}}}\end{picture}}%
\end{picture}}}\end{picture}}%

\newcommand{\NBIMO}[1]%
{\begin{picture}(0,0)%
\put(0,0){\makebox(0,0)%
{\begin{picture}(0,#1)%
\put(0,0){\vector(0,1){#1}}%
\truex{300}\truey{600}\truez{800}%
\put(-\value{x},\value{x}){\line(1,0){\value{y}}}%
\put(0,#1){\begin{picture}(0,0)%
\put(-\value{x},-\value{z}){\line(1,0){\value{y}}}\end{picture}}%
\end{picture}}}\end{picture}}%

\newcommand{\NBIAR}[1]%
{\begin{picture}(0,0)%
\truex{200}%
\put(0,0){\makebox(0,0)%
{\begin{picture}(0,#1)\put(-\value{x},0){\vector(0,1){#1}}%
\put(\value{x},0){\vector(0,1){#1}}%
\end{picture}}}\end{picture}}%

\newcommand{\Nisov}[2]%
{\Vbicase{\NAR}{#1\hspace{-2pt}}{\hspace{-2pt}\cong}{#200}}%

\newcommand{\fdcase}[3]{\begin{picture}(0,0)%
\put(0,-150){#1}%
\truex{200}\truey{600}\truez{600}%
\put(-\value{x},-\value{x}){\makebox(0,\value{z})[r]{${#2}$}}%
\put(\value{x},-\value{y}){\makebox(0,\value{z})[l]{${#3}$}}%
\end{picture}}%

\newcommand{\NEAR}{\begin{picture}(0,0)%
\put(-2900,-2900){\vector(1,1){5800}}%
\end{picture}}%

\newcommand{\NEDOTAR}%
{\truex{100}\truey{212}%
\NUMBEROFDOTS=5800%
\divide\NUMBEROFDOTS by \value{y}%
\begin{picture}(0,0)%
\multiput(-2900,-2900)(\value{y},\value{y}){\NUMBEROFDOTS}%
{\circle*{\value{x}}}%
\put(2900,2900){\vector(1,1){0}}%
\end{picture}}%

\newcommand{\near}{\fdcase{\NEAR}{}{}}%

\newcommand{\SWDOTAR}%
{\truex{100}\truey{212}%
\NUMBEROFDOTS=5800%
\divide\NUMBEROFDOTS by \value{y}%
\begin{picture}(0,0)%
\multiput(2900,2900)(-\value{y},-\value{y}){\NUMBEROFDOTS}%
{\circle*{\value{x}}}%
\put(-2900,-2900){\vector(-1,-1){0}}%
\end{picture}}%

\newcommand{\sdcase}[3]{\begin{picture}(0,0)%
\put(0,-150){#1}%
\truex{100}\truez{600}%
\put(\value{x},\value{x}){\makebox(0,\value{z})[l]{${#2}$}}%
\truex{300}\truey{800}%
\put(-\value{x},-\value{y}){\makebox(0,\value{z})[r]{${#3}$}}%
\end{picture}}%

\newcommand{\SEAR}{\begin{picture}(0,0)%
\put(-2900,2900){\vector(1,-1){5800}}%
\end{picture}}%

\newcommand{\SEDOTAR}%
{\truex{100}\truey{212}%
\NUMBEROFDOTS=5800%
\divide\NUMBEROFDOTS by \value{y}%
\begin{picture}(0,0)%
\multiput(-2900,2900)(\value{y},-\value{y}){\NUMBEROFDOTS}%
{\circle*{\value{x}}}%
\put(2900,-2900){\vector(1,-1){0}}%
\end{picture}}%

\newcommand{\sear}{\sdcase{\SEAR}{}{}}%

\newcommand{\NWDOTAR}%
{\truex{100}\truey{212}%
\NUMBEROFDOTS=5800%
\divide\NUMBEROFDOTS by \value{y}%
\begin{picture}(0,0)%
\multiput(2900,-2900)(-\value{y},\value{y}){\NUMBEROFDOTS}%
{\circle*{\value{x}}}%
\put(-2900,2900){\vector(-1,1){0}}%
\end{picture}}%

\newcommand{\ENEAR}[2]%
{\makebox[0pt]{\begin{picture}(0,0)%
\put(0,-150){\makebox(0,0){\begin{picture}(0,0)%
\put(-6600,-3300){\vector(2,1){13200}}%
\truex{200}\truey{800}\truez{600}%
\put(-\value{x},\value{x}){\makebox(0,\value{z})[r]{${#1}$}}%
\put(\value{x},-\value{y}){\makebox(0,\value{z})[l]{${#2}$}}%
\end{picture}}}\end{picture}}}%

\newcommand{\ESEAR}[2]%
{\makebox[0pt]{\begin{picture}(0,0)%
\put(0,-150){\makebox(0,0){\begin{picture}(0,0)%
\put(-6600,3300){\vector(2,-1){13200}}%
\truex{200}\truey{800}\truez{600}%
\put(\value{x},\value{x}){\makebox(0,\value{z})[l]{${#1}$}}%
\put(-\value{x},-\value{y}){\makebox(0,\value{z})[r]{${#2}$}}%
\end{picture}}}\end{picture}}}%

\newcommand{\WNWAR}[2]%
{\makebox[0pt]{\begin{picture}(0,0)%
\put(0,-150){\makebox(0,0){\begin{picture}(0,0)%
\put(6600,-3300){\vector(-2,1){13200}}%
\truex{200}\truey{800}\truez{600}%
\put(\value{x},\value{x}){\makebox(0,\value{z})[l]{${#1}$}}%
\put(-\value{x},-\value{y}){\makebox(0,\value{z})[r]{${#2}$}}%
\end{picture}}}\end{picture}}}%

\newcommand{\WSWAR}[2]%
{\makebox[0pt]{\begin{picture}(0,0)%
\put(0,-150){\makebox(0,0){\begin{picture}(0,0)%
\put(6600,3300){\vector(-2,-1){13200}}%
\truex{200}\truey{800}\truez{600}%
\put(-\value{x},\value{x}){\makebox(0,\value{z})[r]{${#1}$}}%
\put(\value{x},-\value{y}){\makebox(0,\value{z})[l]{${#2}$}}%
\end{picture}}}\end{picture}}}%

\newcommand{\NNEAR}[2]%
{\raisebox{-1pt}[0pt][0pt]{\begin{picture}(0,0)%
\put(0,0){\makebox(0,0){\begin{picture}(0,0)%
\put(-3300,-6600){\vector(1,2){6600}}%
\truex{100}\truez{600}%
\put(-\value{x},\value{x}){\makebox(0,\value{z})[r]{${#1}$}}%
\put(\value{x},-\value{z}){\makebox(0,\value{z})[l]{${#2}$}}%
\end{picture}}}\end{picture}}}%

\newcommand{\SSWAR}[2]%
{\raisebox{-1pt}[0pt][0pt]{\begin{picture}(0,0)%
\put(0,0){\makebox(0,0){\begin{picture}(0,0)%
\put(3300,6600){\vector(-1,-2){6600}}%
\truex{100}\truez{600}%
\put(-\value{x},\value{x}){\makebox(0,\value{z})[r]{${#1}$}}%
\put(\value{x},-\value{z}){\makebox(0,\value{z})[l]{${#2}$}}%
\end{picture}}}\end{picture}}}%

\newcommand{\SSEAR}[2]%
{\raisebox{-1pt}[0pt][0pt]{\begin{picture}(0,0)%
\put(0,0){\makebox(0,0){\begin{picture}(0,0)%
\put(-3300,6600){\vector(1,-2){6600}}%
\truex{200}\truez{600}%
\put(\value{x},\value{x}){\makebox(0,\value{z})[l]{${#1}$}}%
\put(-\value{x},-\value{z}){\makebox(0,\value{z})[r]{${#2}$}}%
\end{picture}}}\end{picture}}}%

\newcommand{\NNWAR}[2]%
{\raisebox{-1pt}[0pt][0pt]{\begin{picture}(0,0)%
\put(0,0){\makebox(0,0){\begin{picture}(0,0)%
\put(3300,-6600){\vector(-1,2){6600}}%
\truex{200}\truez{600}%
\put(\value{x},\value{x}){\makebox(0,\value{z})[l]{${#1}$}}%
\put(-\value{x},-\value{z}){\makebox(0,\value{z})[r]{${#2}$}}%
\end{picture}}}\end{picture}}}%

\newcommand{\Necurve}[2]%
{\begin{picture}(0,0)%
\truex{1300}\truey{2000}\truez{200}%
\put(0,\value{x}){\oval(#200,\value{y})[t]}%
\put(0,\value{x}){\makebox(0,0){\begin{picture}(#200,0)%
\put(#200,0){\vector(0,-1){\value{z}}}%
\put(0,0){\line(0,-1){\value{z}}}\end{picture}}}%
\truex{2500}%
\put(0,\value{x}){\makebox(0,0)[b]{${#1}$}}%
\end{picture}}%

\newcommand{\Nwcurve}[2]%
{\begin{picture}(0,0)%
\truex{1300}\truey{2000}\truez{200}%
\put(0,\value{x}){\oval(#200,\value{y})[t]}%
\put(0,\value{x}){\makebox(0,0){\begin{picture}(#200,0)%
\put(#200,0){\line(0,-1){\value{z}}}%
\put(0,0){\vector(0,-1){\value{z}}}\end{picture}}}%
\truex{2500}%
\put(0,\value{x}){\makebox(0,0)[b]{${#1}$}}%
\end{picture}}%

\newcommand{\Securve}[2]%
{\begin{picture}(0,0)%
\truex{1300}\truey{2000}\truez{200}%
\put(0,-\value{x}){\oval(#200,\value{y})[b]}%
\put(0,-\value{x}){\makebox(0,0){\begin{picture}(#200,0)%
\put(#200,0){\vector(0,1){\value{z}}}%
\put(0,0){\line(0,1){\value{z}}}\end{picture}}}%
\truex{2500}%
\put(0,-\value{x}){\makebox(0,0)[t]{${#1}$}}%
\end{picture}}%

\newcommand{\Swcurve}[2]%
{\begin{picture}(0,0)%
\truex{1300}\truey{2000}\truez{200}%
\put(0,-\value{x}){\oval(#200,\value{y})[b]}%
\put(0,-\value{x}){\makebox(0,0){\begin{picture}(#200,0)%
\put(#200,0){\line(0,1){\value{z}}}%
\put(0,0){\vector(0,1){\value{z}}}\end{picture}}}%
\truex{2500}%
\put(0,-\value{x}){\makebox(0,0)[t]{${#1}$}}%
\end{picture}}%

\newcommand{\Escurve}[2]%
{\begin{picture}(0,0)%
\truex{1400}\truey{2000}\truez{200}%
\put(\value{x},0){\oval(\value{y},#200)[r]}%
\put(\value{x},0){\makebox(0,0){\begin{picture}(0,#200)%
\put(0,0){\vector(-1,0){\value{z}}}%
\put(0,#200){\line(-1,0){\value{z}}}\end{picture}}}%
\truex{2500}%
\put(\value{x},0){\makebox(0,0)[l]{${#1}$}}%
\end{picture}}%

\newcommand{\Encurve}[2]%
{\begin{picture}(0,0)%
\truex{1400}\truey{2000}\truez{200}%
\put(\value{x},0){\oval(\value{y},#200)[r]}%
\put(\value{x},0){\makebox(0,0){\begin{picture}(0,#200)%
\put(0,0){\line(-1,0){\value{z}}}%
\put(0,#200){\vector(-1,0){\value{z}}}\end{picture}}}%
\truex{2500}%
\put(\value{x},0){\makebox(0,0)[l]{${#1}$}}%
\end{picture}}%

\newcommand{\Wscurve}[2]%
{\begin{picture}(0,0)%
\truex{1300}\truey{2000}\truez{200}%
\put(-\value{x},0){\oval(\value{y},#200)[l]}%
\put(-\value{x},0){\makebox(0,0){\begin{picture}(0,#200)%
\put(0,0){\vector(1,0){\value{z}}}%
\put(0,#200){\line(1,0){\value{z}}}\end{picture}}}%
\truex{2400}%
\put(-\value{x},0){\makebox(0,0)[r]{${#1}$}}%
\end{picture}}%

\newcommand{\Wncurve}[2]%
{\begin{picture}(0,0)%
\truex{1300}\truey{2000}\truez{200}%
\put(-\value{x},0){\oval(\value{y},#200)[l]}%
\put(-\value{x},0){\makebox(0,0){\begin{picture}(0,#200)%
\put(0,0){\line(1,0){\value{z}}}%
\put(0,#200){\vector(1,0){\value{z}}}\end{picture}}}%
\truex{2400}%
\put(-\value{x},0){\makebox(0,0)[r]{${#1}$}}%
\end{picture}}%

\newcount\SCALE%

\newcount\NUMBER%

\newcount\LINE%

\newcount\COLUMN%

\newcount\WIDTH%

\newcount\SOURCE%

\newcount\ARROW%

\newcount\TARGET%

\newcount\ARROWLENGTH%

\newcount\NUMBEROFDOTS%

\newcounter{x}%

\newcounter{y}%

\newcounter{z}%

\newcounter{horizontal}%

\newcounter{vertical}%

\newskip\itemlength%

\newskip\firstitem%

\newskip\seconditem%

\newcommand{\printarrow}{}%

\newcommand{\truex}[1]{%
\NUMBER=#1%
\multiply\NUMBER by 100%
\divide\NUMBER by \SCALE%
\setcounter{x}{\NUMBER}}%

\newcommand{\truey}[1]{%
\NUMBER=#1%
\multiply\NUMBER by 100%
\divide\NUMBER by \SCALE%
\setcounter{y}{\NUMBER}}%

\newcommand{\truez}[1]{%
\NUMBER=#1%
\multiply\NUMBER by 100%
\divide\NUMBER by \SCALE%
\setcounter{z}{\NUMBER}}%

\newcommand{\changecounters}[1]{%
\SOURCE=\ARROW%
\ARROW=\TARGET%
\settowidth{\itemlength}{#1}%
\ifdim \itemlength > 2800\unitlength%
\addtolength{\itemlength}{-2800\unitlength}%
\TARGET=\itemlength%
\divide\TARGET by 1310%
\multiply\TARGET by 100%
\divide\TARGET by \SCALE%
\else%
\TARGET=0%
\fi%
\ARROWLENGTH=5000%
\advance\ARROWLENGTH by -\SOURCE%
\advance\ARROWLENGTH by -\TARGET%
\advance\SOURCE by -\TARGET}%

\newcommand{\initialize}[1]{%
\LINE=0%
\COLUMN=0%
\WIDTH=0%
\ARROW=0%
\TARGET=0%
\changecounters{#1}%
\renewcommand{\printarrow}{#1}%
\begin{center}%
\vspace{10pt}%
\begin{picture}(0,0)}%

\newcommand{\DIAGV}[2]{%
\SCALE=#1%
\setlength{\unitlength}{655sp}%
\multiply\unitlength by \SCALE%
\divide\unitlength by 100%
\initialize{\mbox{$#2$}}}%

\newcommand{\n}[1]{%
\changecounters{\mbox{$#1$}}%
\put(\COLUMN,\LINE){\makebox(0,0){\printarrow}}%
\thinlines%
\renewcommand{\printarrow}{\mbox{$#1$}}%
\advance\COLUMN by 4000}%

\newcommand{\nn}[1]{%
\put(\COLUMN,\LINE){\makebox(0,0){\printarrow}}%
\thinlines%
\ifnum \WIDTH < \COLUMN%
\WIDTH=\COLUMN%
\else%
\fi%
\advance\LINE by -4000%
\COLUMN=0%
\ARROW=0%
\TARGET=0%
\changecounters{\mbox{$#1$}}%
\renewcommand{\printarrow}{\mbox{$#1$}}}%

\newcommand{\conclude}{%
\put(\COLUMN,\LINE){\makebox(0,0){\printarrow}}%
\thinlines%
\ifnum \WIDTH < \COLUMN%
\WIDTH=\COLUMN%
\else%
\fi%
\setcounter{horizontal}{\WIDTH}%
\setcounter{vertical}{-\LINE}%
\end{picture}}%

\newcommand{\diag}{%
\conclude%
\raisebox{0pt}[0pt][\value{vertical}\unitlength]{}%
\hspace*{\value{horizontal}\unitlength}%
\vspace{10pt}%
\end{center}%
\setlength{\unitlength}{1pt}}%

\newcommand{\diagv}[3]{%
\conclude%
\NUMBER=#1%
\rule{0pt}{\NUMBER pt}%
\hspace*{-#2pt}%
\raisebox{0pt}[0pt][\value{vertical}\unitlength]{}%
\hspace*{\value{horizontal}\unitlength}
\NUMBER=#3%
\advance\NUMBER by 10%
\vspace*{\NUMBER pt}%
\end{center}%
\setlength{\unitlength}{1pt}}%

\newcommand{\N}[1]%
{\raisebox{0pt}[7pt][0pt]{$#1$}}%

\newcommand{\crosslength}[2]{%
\settowidth{\firstitem}{#1}%
\settowidth{\seconditem}{#2}%
\ifdim\firstitem < \seconditem%
\itemlength=\seconditem%
\else%
\itemlength=\firstitem%
\fi%
\divide\itemlength by 2%
\hspace{\itemlength}}%

\marginsize{2.5cm}{2.5cm}{3cm}{4cm}
\newtheorem{thm}{Theorem}[section]
\newtheorem{cor}[thm]{Corollary}
\newtheorem{lem}[thm]{Lemma}
\newtheorem{prop}[thm]{Proposition}
\newtheorem{ex}{Example}
\theoremstyle{definition}
\newtheorem{defn}[thm]{Definition}
\theoremstyle{remark}

\numberwithin{equation}{section}

\def\Mod{\mbox{-Mod}}

\newcommand{\Natur}{{\mathbb N}}

\def\Ker{\mbox{\rm Ker}}
\def\Hom{\mbox{\rm Hom}}
\def\Sup{\mbox{\rm sup}}
\def\Tor{\mbox{\rm Tor}}
\def\Q{{\mathbb Q}}
\def\Z{{\mathbb Z}}
\def\Ext{\mbox{\rm Ext}}

\def\li{{\displaystyle \lim_{\rightarrow}}\ }
\begin{document}
\baselineskip=28pt
\title[Injective representations of quivers]{Injective representations of infinite quivers.
Applications}%
\address{Department of Mathematics, University of Kentucky, Lexington, Kentucky
40506-0027, U.S.A.} \email{enochs@ms.uky.edu}

\author{E. Enochs, S. Estrada and J.R.Garc\'{\i}a Rozas}
\address{
Departamento de Matem\'atica Aplicada, Universidad de Murcia,Campus
del Espinardo, Espinardo (Murcia) 30100, Spain}
\email{sestrada@um.es}%
\address{
Departamento de \'Algebra y A. Matem\'atico, Universidad de
Almer\'{\i}a, Almer\'{\i}a 04071, Spain}\email{jrgrozas@ual.es}

\keywords{
\\ The authors are partially supported by the DGI MTM2005-03227}%

\thanks{Estrada's work was supported by a MEC/Fulbright grant from the Spanish Secretar\'{\i}a de Estado
de Universidades e Investigaci\'on del Ministerio de Educaci\'on y Ciencia}%
\subjclass{}%

\begin{abstract}
In this article we study injective representations of infinite
quivers. We classify the indecomposable injective representations of
trees and then describe Gorenstein injective and projective
representations of barren trees.
\end{abstract}
\maketitle
\section{Introduction}
The major impetus for the study of the representations of quivers
was given by Gabriel's study of the finite quivers of finite
representation type (see \cite{Gabriel}) and their connection with
the Dynkin diagrams associated with finite dimensional semisimple
Lie algebras over the field of complex numbers.

The classical representation theory of quivers involved finite
quivers and assumed that the ring was an algebraically closed field
with the assumption that all vector spaces involved were finite
dimensional. In \cite{procesi} the study of semisimple
representations of these kinds of quivers were considered. Recently
representations by modules over more general quivers have been
studied. Our concern will be these kinds of representations and is a
continuation of the program initiated in \cite{EnHer} and continued
in \cite{EnLoyPark}, \cite{EnEsGar}, \cite{EnEst} and
\cite{EnLoyBlas}.

Many categories of graded modules over graded rings are equivalent
to representations of quivers which are often infinite. For a simple
example, the category of $\mathbb{Z}$-graded modules over the graded
ring $R[x]=R+R x+R x^2+\cdots $ (here $R$ is any ring with identity)
is equivalent to the category of representations over $R$ of the
infinite line quiver. Less trivial examples can be given involving
group rings $R[G]$ with the obvious grading. In \cite{Estr} it was
shown that the category of quasi-coherent sheaves over any scheme is
equivalent to a category of representations of a quiver (with
certain modifications on the representations). In many of these
cases the quiver viewpoint leads to simplifications of proofs and of
the descriptions of objects in related categories. This has
certainly proved true in the case of finite quivers and promises to
be so when the quivers are infinite.

Our techniques are necessarily different from the usual ones
concerning quivers without oriented cycles, and in general for the
classic treatment of representation theory of associative algebras
(see for example \cite{assem}). We consider the geometric properties
of the quiver and also have used infinite matrix techniques to
classify projective representations (see \cite{EnEst}). In
\cite{EnLoyBlas} flat representations were studied and in
\cite{EnLoyPark} the so called noetherian quivers were
characterized.

So our main concern in this paper is the study of injective
representations over a possibly infinite quiver in terms of local
properties of the representations (see properties $i)$ and $ii)$ of
Proposition \ref{epies}). As we will see in Section
\ref{applications} these local properties turn to be very useful in
studying and characterizing Gorenstein injective, projective and
flat representations of quivers. These are of particular interest
for defining a version of relative homological algebra that is
called Gorenstein homological algebra. For it has been recently
proved in \cite{EnEsGar} (by using the results of \cite{gorenstein})
that a fruitful version of Gorenstein homological algebra can be
developed in the category of representations over an arbitrary
quiver, when the base ring is Gorenstein. Moreover by
\cite{gorenstein} we get that two model structures can be derived in
the category of representations over a quiver where we use these
Gorenstein injective and Gorenstein projective representations to
define the cofibrations (see \cite{hovey} for a deep study on model
category structures).

On the other hand, infinite and barren trees (see Section
\ref{seced}) appear naturally in the study of finiteness conditions
in the category of representations of quivers. In particular
concerning the question of when this category admits a family of
noetherian generators (cf. \cite{EnLoyPark}) and in the
characterization of a Gorenstein path ring (cf. \cite{EnEsGar}). So
in Section \ref{indesec} we study indecomposable injective
representations of trees to determinate the structure of injective
representations on barren trees (see the remark after Corollary
\ref{inybarren}).

Sections \ref{rootedder} and \ref{seced} are devoted to exhibiting a
wide class of quivers whose injective representations admit the
previous ``local'' characterization (the so called source injective
representation quivers, see Definiton \ref{miles}). One interesting
question is to ask if this characterization carries over to quivers
with relations. For instance, in \cite{Estrad} it is proved that
this is the case when we consider the infinite line quiver on both
sides and then take the monomial relation given by the composition
of $N$-paths. So we get a local structure theorem for injective
representations in the category of $N$-complexes of modules ($N\geq
2)$) and all machinery of Section \ref{applications} can be applied
in this situation.

\section{Preliminaries}

Throughout this article we will use the terminology and results of
\cite{EnHer}.

We keep the notation introduced in \cite{assem}. All rings considered
in this paper will be associative with identity and, unless
otherwise specified, not necessarily commutative. The letter $R$
will usually denote a ring.

As usual we denote a quiver by $Q$ with the understanding that a
quiver is a fourtuple $Q=(V(Q),\Gamma(Q),s,t)$ where $V(Q)$ is a set
of vertices, a set $\Gamma(Q)$ of arrows between these vertices and
two maps $s,t:\Gamma(Q)\to V(Q)$, where for each $a\in \Gamma(Q)$,
$s(a)$ and $t(a)$ assign to an arrow $a$ the source vertex and
terminal vertex of $a$ respectively. All quivers considered in this
paper may be infinite, that is, one of the two sets $V(Q)$ or
$\Gamma(Q)$ can be infinite. Note that we do not exclude loops or
multiple arrows in the definition of a quiver. Sometimes we will
denote $V(Q)$ (resp. $\Gamma(Q)$) simply as $V$ (resp. $\Gamma$) if
the quiver is understood. A finite path $p$ of a quiver $A$ is a
sequence of arrows $a_n\cdots a_2a_1$ with $t(a_i)=s(a_{i+1})$ for
all $i=0,\ldots, n-1$. Therefore $s(p)=s(a_1)$ and $t(p)=t(a_n)$.
Two paths $p$ and $q$ of a quiver can be composed, getting another
path $qp$ whenever $t(p)=s(q)$. A quiver can be thought as a small
category where the objects are the vertices and the morphisms are
the paths. The vertices of $Q$ can be considered as the identities
of $Q$, that is, a vertex $v$ of $Q$ is a trivial path where
$s(v)=t(v)=v$. A tree is a quiver $T$ having a vertex $v$ such that
for another vertex $w$ of $T$ there exists a unique path $p$ such
that $s(p)=v$ and $t(p)=v$. Such vertex $v$ is called the root of
the tree $T$. A forest is a quiver in which every connected
component is a tree. For a path $p$ of $Q$ we denote by $t(p)$
(resp. $s(p)$) the final (resp. the initial) vertex of $p$.

A representation by modules $\mathcal{X}$ of a given quiver $Q$ is a
functor ${\mathcal{X}}:Q\to R\Mod$. Such a representation is
determined by giving a module ${\mathcal{X}}(v)$ to each vertex $v$
of $Q$ and a homomorphism ${\mathcal{X}}(a):{\mathcal{X}}(v_1)\to
{\mathcal{X}}(v_2)$ to each arrow $a:v_1\to v_2$ of $Q$. A morphism
$\eta$ between two representations $\mathcal{X}$ and $\mathcal{Y}$
is a natural transformation, so it will be a family
$\{\eta_v\}_{v\in V}$ such that ${\mathcal{Y}}(a)\circ
\eta_{v_1}=\eta_{v_2}\circ {\mathcal{X}}(a)$ for any arrow $a:v_1\to
v_2$ of $Q$. Thus the representations of a quiver $Q$ by modules
over a ring $R$ is a category, denoted by $(Q,R\Mod)$, which is a
Grothendieck category with enough projectives.

The category $(Q,R\Mod)$ is equivalent to the category of modules
over the path algebra $RQ$, which is a ring with enough idempotents
that in general does not have a unit (unless $|V|$ is finite).

For a given quiver one can find a family injective cogenerators from
an adjoint situation as it is shown in \cite{EnHer}. For every
vertex $v\in V$ and the embedding morphism $\{v\}\subseteq Q$ the
family $\{e_*^v(E):\ v\in V\}$ is a family of injective cogenerators
of $(Q,R\Mod)$, whenever $E$ is an injective cogenerator of $R\Mod$.
The functor $e_*^v:R\Mod\to (Q,R\Mod)$ is defined in \cite[Section
4]{EnHer} as $e_*^v(M)(w)=\prod_{Q(w,v)}M$, where $Q(w,v)$ denotes
the set of paths $p$ in $Q$ such that $s(p)=w$ and $t(p)=v$. If
$a:w_1\to w_2$ is an arrow, then $$e_*^v(M)(a):\prod_{Q(w_1,v)}M\to
\prod_{Q(w_2,v)}M $$ is given by the coordinate-wise function
$e_*^v(M)(a)=\prod_{Q(w_2,v)a}id_{Q(w_2,v)}$. Then by \cite[Theorem
4.1]{EnHer} $e_*^v$ is the right adjoint functor of the evaluation
functor $T_v:(Q,R\Mod)\to R\Mod$ given by
$T_v({\mathcal{X}})=\mathcal{X}(v)$ for any representation
$\mathcal{X}\in (Q,R\Mod)$.

We will need the property satisfied by injective representations
given by the next result.
\begin{prop}\label{epies}
Let $Q$ be any quiver and let $(Q,R\Mod)$ the category of
representations of $Q$ by left $R$-modules. If ${\mathcal{X}}\in
(Q,R\Mod)$ is injective then:
\begin{enumerate}
\item[i) ] ${\mathcal{X}}(v)$ is an injective $R$-module, for any vertex $v$ of
$Q$. \item[ii) ] For any vertex $v$ the morphism
$${\mathcal{X}}(v)\to \prod_{s(a)=v}{\mathcal{X}}(t(a))$$ induced by
${\mathcal{X}}(v)\to {\mathcal{X}}(t(a))$ is a splitting
epimorphism.
\end{enumerate}
\end{prop}
{\bf Proof.} We consider the injective representations $e_*^v(E)$
associated with a given vertex of $Q$ and an injective left
$R$-module $E$. By the construction given in the proof of
\cite[Theorem 4.1]{EnHer} we see that each $e_*(E)$ has the property
of the Proposition. As noted in \cite[pg.303]{EnHer} the $e_*^v(E)$
(varying $v$ and $E$) cogenerate $(Q,R\Mod)$. So any injective
representation is a retract of products of the various $e_*^v(E)$'s.
The property is preserved under taking products and retracts and so
we get the desired result. {\hfill$\Box$}

\bigskip\par
Since the category of representations of a quiver $Q$ is always a
Grothendieck category, it has enough injective representations. In
this paper we are mainly concerned with studying when injective
representations are characterized in terms of conditions $i)$ and
$ii)$ of Proposition \ref{epies}. These conditions involve local
properties of the representation on the vertices and their
corresponding sources, so they motivate the following definition
which is pivotal for the rest of the paper.
\begin{defn}\label{miles}
We say that a quiver $Q$ is a source injective representation quiver
if for any $R$ the injective representations of $(Q,R\Mod)$ can be
characterized in terms of conditions $i)$ and $ii)$ of Proposition
\ref{epies}. We will denote by $\mathfrak{I}$ the class of all
source injective representation quivers.
\end{defn}

One of the main goal of the paper is to study the class
$\mathfrak{I}$ and more precisely to determine which trees belong to
it. For this purpose, it is convenient to start Section
\ref{indesec} by determining the indecomposable injective
representations of trees. As we will see in Section
\ref{applications}, there are important properties which can be
derived from the fact that a quiver is in $\mathfrak{I}$.

Thus in Section \ref{rootedder} we prove that a wide class of
quivers, the so called right rooted quivers (that is those which do
not contain a path of the form $\bullet\to \bullet\to \cdots$) are
in $\mathfrak{I}$. But there are also important non-right rooted
quivers in $\mathfrak{I}$, concretely, in Section \ref{seced} we
prove that $A_{\infty}^{\infty}$, $A_{\infty}^+$ and, more generally
any barren tree are source injective representation quivers.

Of course there is a dual of Definition \ref{miles} by dualizing
properties $i)$ and $ii)$ of Proposition \ref{epies}.

\begin{defn}(dual to Definition \ref{miles})
We say that a quiver $Q$ is a sink projective representation quiver
if for any $R$ the projective representation $\mathcal{X}$ of
$(Q,R\Mod)$ can be characterized in terms of the dual of the
conditions $i)$ and $ii)$ of Proposition \ref{epies}, that is,
\begin{enumerate}
\item[i') ] For each vertex $v\in V$, the module ${\mathcal{X}}(v)$ is projective.

\item[ii') ] For a vertex $v\in V$, the morphism
$\oplus_{t(a)=v}{\mathcal{X}}(s(a))\rightarrow {\mathcal{X}}(v)$
(where ${\mathcal{X}}(s(a))\rightarrow {\mathcal{X}}(v)$ is
${\mathcal{X}}(a)$) is a splitting monomorphism.
\end{enumerate}
We will denote by $\mathfrak{P}$ the subclass of the class of all
sink projective representation quivers.
\end{defn}
Left rooted quivers are examples of quivers in $\mathfrak{P}$, but
as it is shown in \cite[Theorem 4.1]{EnEst}, $A_{\infty}^{\infty}$
is not a sink projective representation quiver. Cyclic quivers are
also examples of quivers which do not lie in $\mathfrak{P}$.

\section{Indecomposable injective representations.}\label{indesec}

The main aim of this section is to characterize the indecomposable
injective representations of trees. So let $T$ be a tree quiver. We
recall that for a given $v$ as above, $e_*$ is the right adjoint of
the restriction functor $(T, R\Mod)\to (\{v\},R\Mod)$ (where the
last category is essentially $R\Mod$).

Since our quiver is a tree, these representations have an especially
simple form. These are such that $e^v_*(E)(w)=E$ if there is a path
(necessarily unique) from $w$ to $v$ and such that  $e^v_*(E)(w)=0$
otherwise. And $e^v_*(E)(a)=id_E$ for any arrow $a$ such that
$e^v_*(E)(t(a))=e^v_*(E)(s(a))=E$. This follows from Theorem 4.1 of
\cite{EnHer}. It is easy to see that if $E$ is an indecomposable
injective module then each $e_*^v(E)$ is also an indecomposable
injective object of $(T,R\Mod)$. But in general these $e_*^v(E)$ do
not give all such objects. However, if we introduce the notion of
vertices $w$ at infinity of $T$ and modify the construction of
$e_*^v(E)$ to include objects $e_*^w(E)$ for such vertices at
infinity, we will get all the indecomposable injective objects.

We begin by noting that given the vertex $v$ of $T$ there is a
unique path $p=a_n a_{n-1}\cdots a_1$ of $Q$ such that $s(a_1)$ is
the root of $T$ and such that $t(a_n)=v$. Then $e_*^v(E)$ is such
that $e_*^v(E)(a_i)=id_E$ for $i=1,2,\cdots, n$ and such that
$e_*^v(v')=0$ when $v'$ is none of $s(a_1), t(a_1), \cdots,
s(a_n),t(a_n)$. We consider infinite paths $p=\cdots a_3a_2a_1$ such
that $s(a_1)$ is the root of $T$. For each such path we associate a
$w$ (distinct from all the vertices of $T$) that we call a vertex at
infinity. With two distinct such paths we associate distinct
vertices at infinity (we note the analogy with the procedure of
adding points at infinity to form a projective space from an affine
space).

Now let $p=\cdots a_3 a_2 a_1$ be an infinite path with $s(a_1)$ the
root of the tree and with $w$ its associated vertex at $\infty$. For
an injective left $R$-module $E$ we define the object $e_*^w(E)$ of
$(T,R\Mod)$ to be that unique object such that $e_*^w(E)(a_i)=id_E$
for each $a_1,a_2,a_3,\cdots$ and such that $e_*^w(E)(v)=0$ if $v$
is any vertex distinct from all the vertices $s(a_j),t(a_j)$ for
$j=1,2,3,\cdots$.

We now shall prove
\begin{prop}\label{referencia1}
If $w$ is a vertex at infinity for $T$ and if $E$ is an injective
left $R$-module then $e_*^w(E)$ is an injective object of
$(T,R\Mod)$.
\end{prop}
{\bf Proof.} We first prove the claim when our tree is
$$A_{\infty}^+\equiv \bullet\to \bullet\to\bullet\to\cdots$$
(the infinite line to the right). Note that here we have one vertex
at infinity, say $w$, and that then for $E$, $e_*^w(E)$ is just
$$E\stackrel{id}{\to}E\stackrel{id}{\to}E\stackrel{id}{\to}E\to\cdots$$
Then note that if $${\mathcal{M}}= M_0\to M_1\to M_2\to \cdots$$is
any object of $(A_{\infty}^+,R\Mod)$ then the morphisms $M\to
e_*^w(E)$ are in one-to-one correspondence with maps $\li M_n\to E$.
So then if ${\mathcal{S}}= S_0\to S_1\to S_2\to \cdots$ is a
subobject of $M$, any morphism $S\to e_*^w(E)$ gives $\li S_n\to E$.
But $E$ is injective and $\li$ is left exact, so $\li S_n\to E$ can
be extended to $\li M_n\to E$. This in turn gives an extension $M\to
e_*^w(E)$ and shows that $e_*^w(E)$ is injective.

Now we consider an arbitrary tree $T$. Given a vertex $w$ at
infinity of $T$ there is an embedding $A_{\infty}^+\subseteq T$ that
identifies $A_{\infty}^+$ with the infinite path $p$ from the root
of $T$ and having $w$ as its vertex at infinity. We now apply the
comment at the bottom of page 302 of \cite{EnHer}. Letting
${\mathcal{X}}$ be our $e_*^w(E)$ (as an object in
$(A_{\infty}^+,R\Mod)$) and letting $f:A_{\infty}^+\to T$ be the
embedding above we see that $f_*(e_*^w(E))$ is precisely what we
denote $e_*^w(E)$ as an object of $(T,R\Mod)$. Since these right
adjoint functors preserve injectivity (see Section 5, pg. 303 of
\cite{EnHer}) we get that $e_*^w(E) $ is injective. {\hfill$\Box$}

\bigskip\par\noindent
Note that if $E$ is also indecomposable then $e_*^w(E) $ is an
indecomposable injective object of $(T,R\Mod)$.
\begin{thm}\label{indes}
If ${\mathcal{X}}$ is an indecomposable injective object of
$(T,R\Mod)$ where $T$ is a tree, then ${\mathcal{X}}\cong e_*^w(E)$
for some vertex $w$ of $T$ (finite or infinite) and some
indecomposable injective left $R$-module $E$.
\end{thm}
{\bf Proof.} We will use Proposition \ref{epies}. Let $v_0$ be the
root of $T$. If ${\mathcal{X}}(v_0)\to
\prod_{s(a)=v_0}{\mathcal{X}}(t(a))$ is not an isomorphism, let $E$
be its kernel. Then $E$ is an injective left $R$-module. Also
$e_*^{v_0}(E)$ is a nonzero injective subobject of ${\mathcal{X}}$.
So since ${\mathcal{X}}$ is indecomposable,
${\mathcal{X}}=e_*^{v_0}(E)$.

So now suppose ${\mathcal{X}}(v_0)\to
\prod_{s(a)=v_0}{\mathcal{X}}(t(a))$ is an isomorphism. For any
arrow $b$ with $s(b)=v_0$ we then let ${\mathcal{X}}(v_0)=E'\oplus
E''$ corresponding to the decomposition $\prod
{\mathcal{X}}(t(a))={\mathcal{X}}(t(b))\oplus \prod_{a\neq
b}{\mathcal{X}}(t(a))$. Consider the subobject ${\mathcal{X}}'$ of
${\mathcal{X}}$ such that ${\mathcal{X'}}(v_0)=E'$ and then such
that for $v\neq v_0$ ${\mathcal{X}}(v)=0$ if the unique path from
$v_0$ to $v$ does not go through $b$ and such that
${\mathcal{X'}}(v)={\mathcal{X}}(v)$ if it does.

Then let ${\mathcal{X}}''$ be the subobject of ${\mathcal{X}}$ such
that ${\mathcal{X''}}(v_0)=E''$ and then such that for $v\neq v_0$,
${\mathcal{X''}}(v)=0$ if the path from $v_0$ to $v$ does go through
$b$ and such that ${\mathcal{X''}}(v)={\mathcal{X}}(v)$ otherwise.
Checking the compatibility conditions we see that in fact
${\mathcal{X}}'$ and ${\mathcal{X}}''$ are subobjects and then that
${\mathcal{X}}={\mathcal{X}}'\oplus {\mathcal{X}}''$. So by the
indecomposability of ${\mathcal{X}}$ we get
${\mathcal{X}}={\mathcal{X}}'$ or ${\mathcal{X}}={\mathcal{X}}''$.

Now since ${\mathcal{X}}(v_0)\to \prod {\mathcal{X}}(t(a))$ is an
isomorphism and since ${\mathcal{X}}(v_0)\neq 0$ there is at least
one $b$ such hat ${\mathcal{X}}(t(b))\neq 0$. Choosing this $b$ we
then get that ${\mathcal{X}}'\neq 0$ and so that
${\mathcal{X}}={\mathcal{X}}'$.

So we note that our argument shows that for exactly one arrow $b$
with $s(b)=v_0$ we get that ${\mathcal{X}}(b)$ is an isomorphism and
that ${\mathcal{X}}(a)=0$ for all $a\neq b$ with $s(a)=v_0$. Let
$t(b)=v_1$ for this $b$. Now consider the subtree $T'\subset T$ of
$T$ with root $v_1$ and which contains all paths of $T$ beginning at
$v_1$. If we repeat the arguments above with $T'$ instead of $T$ we
see that if we are then in the first situation we get
${\mathcal{X}}=e_*^{v_1}({\mathcal{X}}(v_0))$. If this does not
occur then the procedure above will give us a $v_2$ and a
corresponding subtree $T''\subset T$. Continuing we see that either
the procedure stops and we get that
${\mathcal{X}}=e_*^{v_n}({\mathcal{X}}(v_0))$ for some vertex $v_n$.
If the procedure does not stop then the vertices
$v_0,v_1,v_2,\cdots$ will determine an infinite path beginning at
$v_0$ and so a vertex $w$ at infinity such that
${\mathcal{X}}=e_*^w({\mathcal{X}}(v_0))$. {\hfill$\Box$}

\bigskip\par\noindent
{\bf Remark.} The $E$ and $w$ such that ${\mathcal{X}}\cong
e_*^{w}(E)$ are uniquely determined up to isomorphisms.

\section{Injective representations of right rooted
quivers}\label{rootedder}

In this section we prove that a wide class of quivers, the so called
right rooted quivers, is contained in the class $\mathfrak{I}$ of
all source injective representation quivers.

Right rooted quivers can be intuitively defined as follows:
\begin{defn}\cite[Proposition 3.6]{EnLoyBlas}
A quiver $Q$ is right rooted (resp. left rooted quiver) if and only
if there exists no path of the form
$\bullet\to\bullet\to\bullet\to\cdots$ in $Q$ (resp. there exists no
path of the form $\cdots\to \bullet\to \bullet\to \bullet $ in $Q$).
\end{defn}

We use the same notation as that introduced in \cite{EnLoyBlas} to
characterize right rooted quivers in terms of certain subsets of the
set of vertices, that is, for a quiver $Q=(V,\Gamma)$ we shall
define by transfinite induction the following subsets of $V$,
$$V_0=\{ v\in V;\ \mbox{there exists no arrow $a$ of $Q$ with }
s(a)=v\}.$$  For a successor ordinal $\alpha$, we define
$$V_{\alpha}=\{ v\in V^{\alpha-1};\ \mbox{there exists no arrow $a$ of
}Q^{\alpha-1}\ \mbox{with }\ s(a)=v \},$$where
$Q^{\alpha-1}=(V^{\alpha-1},\Gamma^{\alpha-1})$ is the subquiver of
$Q$ with
$$V^{\alpha-1}=V\setminus V_{\alpha-1}$$
and $$\Gamma^{\alpha-1}=\Gamma\setminus\{ a\in \Gamma;\ t(a)\in
V_{\alpha-1}\}.$$ For a limit ordinal we define
$$V_{\omega}={\displaystyle \lim_{\rightarrow}}\
V_{\alpha}=\cup_{\alpha <\omega} V_{\alpha}$$ and
$Q^{\omega}=(V^{\omega},\Gamma^{\omega})$ the subquiver of $Q$,
where $V^{\omega}=V\setminus \cup_{\alpha<\omega}V_{\alpha}$ and
$$\Gamma^{\omega}=\Gamma\setminus\{ a\in \Gamma;\ t(a)\in {\displaystyle
\lim_{\rightarrow}}\ V_{\alpha}\}.$$

Then if a quiver is right rooted there will exist an ordinal number
$\lambda$ such that $V=\cup_{\alpha< \lambda} V_{\alpha}$. The
converse is also true by \cite[Proposition 3.6]{EnLoyBlas}.

\begin{thm}\label{injrooted}
Let $Q$ be a right rooted quiver. Then $Q$ is a source injective
representation quiver.
\end{thm}
{\bf Proof.} Let us see that conditions $i)$ and $ii)$ of
Proposition \ref{epies} are sufficient to get an injective
representation. By transfinite induction on the set of vertices
$V_{\alpha}$, we shall construct a family of injective
subrepresentations of $\mathcal{E}$ whose direct product coincides
with $\mathcal{E}$. Let ${\mathcal E}_0=\prod_{v\in
V_0}e_*^v({\mathcal{E}}(v))$. Now, if $\alpha$ is a successor
ordinal and $v\in V_{\alpha+1}$ it follows, by hypothesis,
$${\mathcal{E}}(v)=(\prod_{s(a)=v}{\mathcal{E}}(t(a)))\oplus E^{\alpha+1}_v$$(where the
product
 $\prod_{s(a)=v}{\mathcal{E}}(t(a))$ is actually the product of
${\mathcal E}_{\mu}(v)$ with $\mu\leq\alpha$). Then we define
$${\mathcal{E}}_{\alpha+1}=\prod_{v\in V_{\alpha+1}} e_*^v(E^{\alpha+1}_v).$$
Notice that ${\mathcal{E}}_{\alpha+1}(v)=0$ for all $v\in V_{\mu}$,
$\mu<\alpha+1$. Now if $\gamma$ is a limit ordinal, we define
$${\mathcal{E}}_{\gamma}=\prod_{\alpha<\gamma}{\mathcal{E}}_{\alpha}.$$Since the quiver
is right rooted it is easy to see that there will exist an ordinal
number $\lambda$ such that
$${\mathcal{E}}=\prod_{\alpha<\lambda}{\mathcal{E}}_{\alpha}$$ and, since every
${\mathcal{E}}_{\alpha}$ is injective ($e_*$ preserves injective
objects), it follows that ${\mathcal{E}}$ will be also injective.
{\hfill$\Box$}

\begin{ex}\label{inyr}
{\rm Consider the quiver $T\equiv\bullet\to \bullet$, and the
representations ${\mathcal{S}}= S_1\stackrel{\beta}{\to} S_0,\
{\mathcal{X}}= {\mathcal{X}}_1\stackrel{\alpha}{\to}X_0$. If
\DIAGV{60} {S_1} \n{\Ear{\beta}} \n{S_0}\nn
          {\Sar{v_1}}\n{}             \n{\saR{v_0}}\nn
          {X_1}      \n{\Ear{\alpha}} \n{X_0}\diag
is a monomorphism of representations and \DIAGV{60} {S_1}
\n{\Ear{\beta}}  \n{S_0}\nn
          {\Sar{h_1}}\n{}             \n{\saR{h_0}}\nn
          {E_1}      \n{\Ear{f}} \n{E_0}\diag
a morphism of representations with $E_1\stackrel{f}{\to}E_0$ an
injective representation, we shall explicitly describe a morphism of
representations $(t_1,t_0)$ extending $(h_1,h_0)$ via $(v_1,v_0)$.
Since $E_0$ is injective there exists $t_0:X_0\to E_0$ such that
$t_0\circ v_0=h_0$. On the other hand, since $f$ splits, let $g$ be
a section. Then it follows that $E_1={\rm Ker}f\oplus {\rm Im} g$,
so we may define a morphism $\gamma:X_1\to {\rm Ker}f$ such that
$\gamma\circ v_1=p_{{\rm Ker}f}\circ h_1$ (where $p_{{\rm Ker}f}$
denotes the canonical projection of $E_1$ over ${\rm Ker}f$ ). Now
we define $t_1:X_1\to E_1$ as $t_1=\gamma+g\circ t_0\circ \alpha$.
It is immediate that $f\circ t_1=t_0\circ \alpha$. Let us see that
$t_1\circ v_1=h_1$, so let $s_1\in S_1$ and suppose $h_1(s_1)=a+b,$
with $a\in {\rm Ker}f$ and $b\in {\rm Im} g$ (so $b=g(x)$ for some
$x\in E_0$). Then $(t_1\circ v_1)(s_1)=(\gamma v_1)(s_1)+(gt_0\alpha
v_1)(s_1)=p_{{\rm
Ker}f}(h_1(s_1))+(gt_0v_0\beta)(s_1)=a+(gh_0\beta)(s_1)=a+(gfh_1)(s_1)
=a+(gf)(a+b)=a+(gf)(g(x))=a+g(x)=a+b=h_1(s_1).$ }
\end{ex}
\begin{ex}\label{extensiones}
{\rm The way of constructing extensions of the previous example can
be easily generalized to many other right rooted quivers. For
example, \DIAGV{30}
           {}\n{}\n{} \n{}     \n{\bullet}\n{} \n{}  \n{} \n{} \n{}  \n{}  \n{\bullet}\nn
           {}\n{}\n{} \n{\near}\n{}\n{} \n{} \n{}  \n{}\n{}\n{}\n{} \n{\sear}       \nn
           {Q_1\equiv}\n{}\n{\bullet}     \n{} \n{} \n{} \n{} \n{}\n{} \n{Q_2\equiv}  \n{}  \n{} \n{}   \n{\bullet}         \nn
           {}\n{} \n{} \n{\sear}\n{}\n{} \n{}\n{}\n{}\n{} \n{}  \n{} \n{\near}                 \nn
           {}\n{}\n{} \n{} \n{\bullet}\n{} \n{}  \n{}\n{} \n{}  \n{}   \n{\bullet}
\diag $Q_3\equiv$
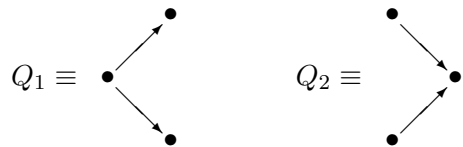
\begin{picture}(40,4)(0,-2)
\put(0,-2){$\bullet$} \put(9,2){\vector(1,0){20}}
 \put(9,-2){\vector(1,0){20}}
 \put(31,-2){$\bullet$}
\end{picture} or more generally to any dual tree or finite
multiple (dual) tree. By the dual tree of a given one we mean a tree
whose arrows are inverted, that is, a quiver with a vertex $v_0$
such that for any other vertex $w$ there exists a unique path from
$w$ to $v_0$. A multiple (dual) tree is defined as a (dual) tree
that admits multiple arrows from a vertex another (so it is not
properly a (dual) tree).}
\end{ex}
\begin{ex}\label{ciclo}
{\rm We now give an example of a non right rooted quiver such that
conditions $i)$ and $ii)$ of Theorem \ref{injrooted} do not imply
that a representation is injective. Let us consider the quiver
 $$Q\equiv 
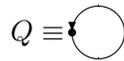
\begin{picture}(20,20)(0,-3)
\put(0,0){\circle*{3}} \put(10,0){\oval(20,20)}
\put(0,1){\vector(0,-1){1}}\end{picture}$$ and the category of
representations of $Q$ over $k$-vector spaces, $(Q,k\Mod)$ ($k$ is a
field). It is easy to see that $(Q,k\Mod)$ is equivalent to the
category $k[x]$-Mod. So now let us consider the representation with
$k[x,x^{-1}]$ in the vertex and the morphism
$k[x,x^{-1}]\stackrel{\cdot x}{\longrightarrow}k[x,x^{-1}]$. It is
obvious that it satisfies $i)$ and $ii)$ of Theorem \ref{injrooted},
but $k[x,x^{-1}]$ is not divisible as a $k[x]$-module (in fact is
$x$-divisible) so can not be injective. }

\end{ex}

\section{ Injective representations of
$A_{\infty}^+$}\label{seced}

In this section we develop new techniques (necessarily very
different from those of Section \ref{rootedder}) to characterize
injective representations of non rooted right quivers. We focus our
attention on the quiver $A_{\infty}^+$ and then generalize this
argument to any quiver such that the connected components are barren
trees (so a forest of barren trees). Barren trees are important
because they appear naturally in the study the existence of
injective covers in $(Q,R\Mod)$ (see \cite[Theorem 3.6]{EnLoyPark})
and also appear in the study of Gorenstein path rings, that is path
rings which are left and right noetherian and whose projective
objects have finite injective dimension (see \cite{iwanaga}). We
recall from \cite{EnLoyPark} the definition of a barren tree: if $T$
is a tree with root $v$, we can divide the set of vertices into
``states" in such a way that the first state $V_1=\{v\}$ and, for
$i\in \Natur,$
$$V_{i+1}=\{w\in V:\ {\rm there\ exists}\ a:v_i\to w,\ \ v_i\in V_i\}.$$ Then
$T$ is barren if the sequence $(n_i)_{i\in \Natur}$ stabilizes
(where $n_i=|V_i|$).

{\bf Example.} The following is an example of an infinite barren
tree.
\bigskip
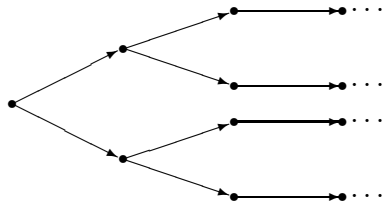
\begin{center}
\begin{picture}(86,86)(0,-36)
\put(0,0){\circle*{3}} \put(0,0){\vector(2,1){40}}
\put(0,0){\vector(2,-1){40}} \put(42,21){\circle*{3}}
\put(42,-21){\circle*{3}} \put(42,21){\vector(3,1){40}}
\put(84,35.3){\circle*{3}}\put(84,35.3){\vector(1,0){40}}\put(125,35.3){\circle*{3}}
\put(128,35.3){\ldots} \put(42,21){\vector(3,-1){40}}
\put(84,6.7){\circle*{3}}\put(84,6.7){\vector(1,0){40}}\put(125,6.7){\circle*{3}}
\put(128,6.7){\ldots} \put(42,-21){\vector(3,1){40}}
\put(84,-6.7){\circle*{3}}\put(84,-6.7){\vector(1,0){40}}\put(125,-6.7){\circle*{3}}
\put(128,-6.7){\ldots} \put(42,-21){\vector(3,-1){40}}
\put(84,-35.3){\circle*{3}}\put(84,-35.3){\vector(1,0){40}}\put(125,-35.3){\circle*{3}}
\put(128,-35.3){\ldots}
\end{picture}
\end{center}

Let $(N_i)_{i\geq 0}$ be any family of modules. Let
$U_n=\prod_{i=n}^{\infty}N_i$ and let ${\mathcal U}=U_0\rightarrow
U_1\rightarrow \cdots$ where each $U_n\rightarrow U_{n+1}$ is the
projection map. Then given ${\mathcal
M}=M_0\stackrel{f_0}{\rightarrow} M_1\stackrel{f_1}{ \rightarrow}
\cdots$ and linear maps $\sigma_n :M_n\rightarrow N_n $ for each
$n\geq 0$ there is a unique corresponding ${\mathcal M}\rightarrow
{\mathcal N}$ such that $M_n\rightarrow U_n=\prod_{i=n}^{\infty}N_i$
is such that the composition $M_n\rightarrow \prod_{i=n}^{\infty}N_i
\rightarrow N_m$ for $m\geq n$ is $\sigma_n$ if $m=n$ and is $
\sigma_m\circ f_{m-1}\circ \cdots \circ f_n $ if $m>n$. And in fact
every morphism ${\mathcal M}\rightarrow {\mathcal N}$ is of this
form. So it is easy to see that if each $N_i$ is injective then
${\mathcal U}$ is injective. Also by choosing the $\sigma$'s
injective we can get an embedding into an injective representation.
So every injective representation is a direct summand of
such a ${\mathcal U}$.\\
Hence we consider direct sum decompositions of arbitrary
representations ${\mathcal M}$, say ${\mathcal M}={\mathcal
M}'\oplus {\mathcal M}''$ where ${\mathcal M}=M_0\rightarrow
\cdots$. Such a decomposition is given by decompositions
$M_n=M_n'\oplus M_n''$ for each $n$, but which are compatible with
the maps $M_n\rightarrow M_{n+1}$, i.e. which map $M_n'$ into
$M_{n+1}'$ and $M_n'' $ into $M_{n+1}''$. So this just says
$M_n\rightarrow M_{n+1}$ is homogeneous with respect to the
decompositions. Given such a decomposition we easily see that each
$\Ker(M_n\rightarrow M_m)$ ($m\geq n$) is a homogeneous submodule
(again with respect to the decomposition). This means that $K=(K\cap
M_n')\oplus
(K\cap M_n'')$ where $K$ is this kernel.\\
 If now in fact ${\mathcal M}$ is such that each $M_n\rightarrow M_{n+1}$ is
surjective then letting ${\mathcal M}=M\rightarrow M/K_0\rightarrow
M/K_1\cdots$ with $K_0\subset K_1 \cdots$ it is not hard to see that
to get a decomposition of ${\mathcal M}$ we only need give a direct
sum decomposition $M=M'\oplus M''$ such that each of the $K_i$ is a
homogeneous submodule. Clearly this remark applies to any
$M_0\rightarrow M_1\rightarrow \cdots $ where all the maps are
surjective. This means that to give a decomposition of this
representation we only need give one of $M_0$ such that all the
kernels $M_0\rightarrow M_n$
are homogeneous for every $n\geq 1$.\\
 We now apply these remarks.

\begin{thm}
Let $(E_i)_{i\geq 0}$ be a family of injective modules. Let
$\oplus_{i=0}^{\infty}E_i \subset E$ be an injective envelope. Then
the representation
$$\overline{\mathcal E}=E\rightarrow E/E_0\rightarrow E/(E_0\oplus
E_1)\rightarrow \cdots$$ is an injective envelope of
$${\mathcal E}=\oplus_{i=0}^{\infty}E_i\rightarrow \oplus_{i=1}^{\infty}E_i\rightarrow
\cdots.$$
\end{thm}
{\bf Proof.} We use the earlier remarks to construct a known
injective representation which is isomorphic to $\overline{\mathcal
E}$. We let ${\mathcal U}$ be the representation $\mathcal U=$
$\prod_{i=0}^{\infty}E_i \rightarrow \prod_{i=1}^{\infty}E_i
\rightarrow \cdots$. This representation is injective. We will try
to find a direct summand which is isomorphic to the representation
of the theorem. So we construct a direct sum decomposition of
${\mathcal U}$. This representation is such that the maps
$\prod_{i=n}^{\infty}E_i\rightarrow \prod_{i=n+1}^{\infty}E_i$ are
all surjective. So we use this fact to construct our decomposition.
Let $\oplus_{i=0}^{\infty}E_i\subset E'\subset
\prod_{i=0}^{\infty}E_i$ where $E'$ is
an injective envelope of $\oplus_{i=0}^{\infty}E_i$.\\
Since $E'$ is injective we have a decomposition
$\prod_{i=0}^{\infty}E_i= E'\oplus E''$. We claim $ker
(\prod_{i=0}^{\infty}E_i \rightarrow \prod_{i=n}^{
\infty}E_i)=E_0\oplus \cdots \oplus E_{n-1}$ is homogeneous with
respect to this decomposition. This is trivial since $E_0\oplus
\cdots E_{n-1}\subset E'$ for each $n\geq 1$. So we now use our
procedure and construct the corresponding decomposition ${\mathcal
E}={\mathcal E}'\oplus {\mathcal E}''$. By
construction it is clear that ${\mathcal E}'$ is isomorphic to the representation\\
$\overline{\mathcal E}=E\rightarrow E/E_0\rightarrow E/(E_0\oplus
E_1) \rightarrow \cdots$. Hence
this representation is injective.\\
To see that it is an envelope, note that $E/(E_0 \oplus \cdots
\oplus E_{n-1})$ is an injective envelope of
$\oplus_{i=n}^{\infty}E_i $ for each $n$. Hence any
$\overline{\mathcal E}\rightarrow \overline{\mathcal E}$ that makes
the obvious diagram commutative gives automorphisms of each
$E/(E_0\oplus \cdots\oplus E_{n-1})$and so ${\mathcal E}\rightarrow
{\mathcal E}$ is an automorphism. {\hfill$\Box$}

\bigskip\par\noindent
   We now want to give a characterization of injective
representations. We will need
\begin{lem}
If $E$ is an injective module then ${\mathcal E}=E\rightarrow
E\rightarrow \cdots $ is an injective representation.
\end{lem}
{\bf Proof.} This follows from the observation (made in the first
part of the proof of Proposition \ref{referencia1}) that morphisms
${\mathcal M}\rightarrow {\mathcal E}$ are given by  maps $\li
M_n\rightarrow E$. {\hfill$\Box$}

\bigskip\par
So we note that any ${\mathcal E}=E_0\rightarrow E_1\rightarrow
\cdots$ with all $E_n$ injective and each $E_n\rightarrow E_{n+1}$
an isomorphism is an
injective representation.\\
We note that there is a natural torsion theory on the category of
representations where ${\mathcal M}=M_0\rightarrow \cdots $ is
torsion if for each $n$ and each $x\in M_n$ there is an $m>n$ such
that $x\in ker(M_n\rightarrow M_m)$. Then for an arbitrary
representation ${\mathcal M}$ we let $t({\mathcal M})$ be the
torsion subrepresentation. Note that $N_0\rightarrow N_1\rightarrow
\cdots$ torsion free just means that each $N_n\rightarrow N_{n+1}$
is an injection.
\begin{thm}\label{deed}
A representation ${\mathcal G}=G_0\rightarrow \cdots $ is injective
if and only if each $G_n$ is injective and each $G_n\rightarrow
G_{n+1}$ is surjective with an injective kernel, that is,
$A_{\infty}^+\in \mathfrak{I}$.
\end{thm}
{\bf Proof.} By Proposition \ref{epies} we know the conditions are
necessary, so assume these conditions. We want to argue that
${\mathcal G}$ is injective. Let $E_0=ker(G_0\rightarrow G_1)$. Then
$E_0$ is injective and is a submodule of the injective
$\Ker(G_0\rightarrow G_2)$. So we can find an $E_1$ so that
$E_0\oplus E_1=\Ker(G_0\rightarrow G_2)$. We proceed in this manner
and find $E_n$ for all $n$ so that $\oplus_{i=0}^{
n-1}E_i=\Ker(G_0\rightarrow G_n)$. Let $\oplus_{i=0}^{\infty}E_i
\subset E \subset G_0$ where $E$ is an injective envelope of
$\oplus_{i=0}^{\infty}E_i$. Then
 up to isomorphism, $\overline{\mathcal E}=E\rightarrow E/E_0\rightarrow
E/(E_0\oplus E_1) \rightarrow \cdots $ is a subrepresentation of
${\mathcal G}$. By the above it is an injective representation and
in fact is the injective envelope of ${\mathcal
E}=\oplus_{i=0}^{\infty}E_i\rightarrow \oplus_{i=1}^{\infty}
E_i\rightarrow \cdots $. Clearly $\mathcal E$ is $t({\mathcal G})$.
So ${\mathcal G}= \overline{\mathcal E}\oplus {\mathcal G}'$. Since
$\overline{\mathcal E} \supset t({\mathcal G})$, ${\mathcal G}'$ is
torsion free. But then ${\mathcal G}'$ has all its terms injective
and all its maps surjective. But this means all its maps are
isomorphisms. Hence ${\mathcal G}'$ is injective and so ${\mathcal
G} $ as the direct sum of two injectives is injective.{\hfill$\Box$}

\bigskip\par
{\bf Remark 1.} From the arguments it is easy to see that an
injective representation is uniquely determined up to isomorphism by
the family $(E_i)_{i\geq 0}$ as above  and the single injective
module $E$ so that the torsion free quotient of the
injective is isomorphic to $E\rightarrow E\rightarrow \cdots $.\\

{\bf Remark 2.} Since a torsion theory is stable if and only if
injective envelopes of torsion objects are torsion we see that if
our theory is stable, then using the notation above,
$\oplus_{i=0}^{\infty}E_i$ must be its own injective envelope for
any family $(E_i)_{i\geq 0}$ of injective modules. So this means the
ring $R$ must be left noetherian. If conversely $R$ is left
noetherian then by our arguments we see that for any injective
${\mathcal G}$, $t({\mathcal G})$ is also injective. This quickly
gives that injective envelopes of torsion objects are also torsion.
So the torsion theory is stable. Hence we get that $R$ is left
noetherian
if and only if our torsion theory is stable.\\

Theorem \ref{deed} allows us to prove that some non right rooted
quivers are also source injective representation quivers.
\begin{cor}\label{inylinea}
Let ${\mathcal{E}}=\cdots\to E_{-2}\to E_{-1}\to E_0\to E_1\to
E_2\to \cdots$ be a representation of the quiver
$$A_{\infty}^{\infty}\equiv\cdots \to \bullet\to \bullet\to \cdots.$$ Then
${\mathcal{E}}$ is injective if, and only if, $E_i$ is an injective
$R$-module and $E_i\to E_{i+1}$ is an splitting epimorphism for all
$ i\in \Z$, that is, the quiver $A_{\infty}^{\infty}$ lies in the
class $\mathfrak{I}$.
\end{cor}
{\bf Proof.} Let $0\to {\mathcal{S}}\stackrel{g}{\to} {\mathcal{X}}$
and ${\mathcal{S}}\stackrel{h}{\to} {\mathcal{E}}$ be a morphism of
representations of $A_{\infty}$ (so $g=(g_i)_{i\in \Z}$ and
$h=(h_i)_{i\in \Z}$). By Theorem \ref{deed} there exist morphisms
$(t_i)_{i\geq 0}$ such that $(h_i)=(t_i)\circ (g_i)$ and $(t_i)$
satisfying the usual compatibility conditions, for all $i\geq 0$.
Now from $t_0:X_0\to E_0$ we get, by Example \ref{inyr}, a morphism
$t_{-1}:X_{-1}\to E_{-1}$ such that $(t_{-1},t_0)$ extends
$(h_{-1},h_0)$ and verifies the compatibility condition. Now proceed
inductively to get a family $t_{-i}:X_{-i}\to E_{-i}$, for all
$i\geq 1$ by applying Example \ref{inyr} iteratively. {\hfill$\Box$}

\bigskip\par
We finish by characterizing injective representations for infinite
barren trees.

\begin{cor} \label{inybarren}
Let $Q$ be a forest whose connected components are barren trees.
Then $Q\in \mathfrak{I}$.
\end{cor}
{\bf Proof.} We need to check that conditions $i)$ and $ii)$ of
Proposition \ref{epies} are sufficient to get an injective
representation.

It is clear that we only need to prove the result for a barren tree
$T$. Since $T$ is barren there exists $k\in \Natur$ such that
$n_k=n_{k+i}$, for all $i\geq 0$. Now suppose that
$V_k=\{v_1,\cdots, v_m,v_{m+1},\cdots,v_{n_k}\}$, where we denote by
$v_1,\cdots, v_m$ the vertices at infinity, so we have
$$v_j\equiv v_j^0\to v_j^1\to v_j^2 \cdots$$ for $1\leq j\leq m$
and suppose (without loss of generality) that $v_j^0\in V_k$, for
all $1\leq j\leq m$.

Now let ${\mathcal{E}}$ be a representation of $T$ satisfying $i)$
and $ii)$ of Proposition \ref{epies} and $0\to
{\mathcal{S}}\stackrel{g}{\to} {\mathcal{X}}$,
${\mathcal{S}}\stackrel{h}{\to}{\mathcal{E}}$ two morphisms of
representations over $T$. Let us consider the corresponding
restrictions $g(v_j)$ and $h(v_j)$ to the vertices at infinity of
$V_k$, $1\leq j\leq m$. Then by Theorem \ref{deed}, for each $1\leq
j\leq m$ there exists a family of morphisms
$$t(v_j)=\{t(v_j^i):\ \ i\in\Natur\}$$ such that $t(v_j^i)\circ
g(v_j^i)=h(v_j^i)$, for all $1\leq j\leq m$, $i\in \Natur$ and for a
fixed $1\leq j\leq m$ the morphisms $\{t(v_j^i):\ \ i\in\Natur\}$
induce a morphism of representations over the quiver $v_j\equiv
v_j^0\to v_j^1\to v_j^2 \cdots$.

Now from the morphisms $$\{t(v_j^0):\ \ 1\leq j\leq m\}\cup
\{t(v_j):\ \ m\leq j\leq n_k\}$$ related to vertices of $V_k$, we
may construct the corresponding extensions to vertices of $V_{k-1}$
as it is shown in Example \ref{inyr} (and by the comments made in
Example \ref{extensiones}) and following repeatedly until we reach
to the root $V_1=\{v\}$. So at the end we will have a representation
$t:{\mathcal{X}}\to {\mathcal{E}}$ extending $h:{\mathcal{S}}\to
{\mathcal{E}}$ via $g:{\mathcal{S}}\to {\mathcal{X}}$ and satisfying
the compatibility condition, that is, if $a:v\to w$ is an arrow of
$T$ the diagram \DIAGV{75} {{\mathcal{X}}(v)}
\n{\Ear{{\mathcal{X}}(a)}} \n{{\mathcal{X}}(w)}\nn
          {\Sar{t(v)}}  \n{}           \n{\saR{t(w)}}\nn
          {{\mathcal{E}}(v)}      \n{\Ear{{\mathcal{E}}(a)}} \n{{\mathcal{E}}(w)}\diag
is commutative.{\hfill$\Box$}

\bigskip\par
{\bf Remark.} If the ring $0\neq R$ is left noetherian and again our
quiver is a tree $T$, then in \cite{EnLoyPark} it was proved that
$(T,R\Mod)$ is a locally noetherian Grothendieck category if and
only if $T$ is barren. In this case every injective object of
$(T,R\Mod)$ is uniquely up to isomorphism the direct sum of the
indecomposable injectives $e_*^w(E)$ of Section \ref{indesec}.

\section{Applications: Gorenstein
representations}\label{applications}


As an application of the results of the previous sections, we study
and characterize the representations of finite injective dimension
(resp. finite projective dimension) and its right orthogonal class
(resp. its left orthogonal class) with respect to the $\Ext^1$
functor, that is, the class of Gorenstein injective (resp.
Gorenstein projective) representations. We also study Gorenstein
flat representations and give and upper bound of the Gorenstein
injective dimension of a representation having finite Gorenstein
injective dimension on every module associated with each vertex.

We recall (see, for example, \cite{EdO}):
\begin{defn}
An object $M$ of an abelian category ${\mathcal A}$ is said to be
Gorenstein injective if there is an exact sequence $$\cdots \to
E_{-3}\to E_{-2}\to E_{-1}\to E_0\to E_{1}\to E_{2}\to \cdots$$ of
injective objects such that $M=\Ker(E_0\to E_{1})$ and such that the
sequence is $\Hom(E,-)$-exact for every injective object $E$ (i.e.
the functor $\Hom(E,-)$ leaves the sequence exact).

\end{defn}

We also recall from \cite{gorenstein} the definition of a Gorenstein
category.

\begin{defn}\label{defgorens}
We will say that a Grothendieck category ${\mathcal{A}}$ is a
Gorenstein category if the following hold:

1) For any object $L$ of ${\mathcal{A}}$, $projdim\, L<\infty$ if
and only if $injdim\, L<\infty$.

2) Finitistic projective dimension and finitistic injective
dimension of ${\mathcal{A}}$ are both finite, that is,
 $$FPD({\mathcal{A}})=sup\{projdim(M):\
projdim(M)<\infty\} $$ and
$$FID({\mathcal{A}})=sup\{injdim(M):\ injdim(M)<\infty\}$$

3) ${\mathcal{A}}$ has a generator $L$ such that $projdim\,
L<\infty$.

\end{defn}

Then in \cite{EnEsGar} it is proved that, for any arbitrary quiver
$Q$, if $R$ is such that $R\Mod$ is Gorenstein (for example if $R$
is an $n$-Gorenstein ring) then $(Q,R\Mod)$ is a Gorenstein category
and so by \cite[Theorem 2.25]{gorenstein} the pair of classes
$(\mathfrak{L},\mathfrak{L}^{\perp})$ where $\mathfrak{L}$ is the
class of representations of finite injective dimension and
$\mathfrak{L}^{\perp}$ is the class of Gorenstein injective
representations, is a complete cotorsion pair. Analogously, the pair
$(^{\perp}\mathfrak{D},\mathfrak{D})$ where $\mathfrak{D} $ is the
class of objects of finite projective dimension and
$^{\perp}\mathfrak{D}$ is the class of Gorenstein projective
representations, is also a complete cotorsion pair.

We now characterize the representations of the classes
$\mathfrak{L}$ and $\mathfrak{L}^{\perp}$ (resp.
$^{\perp}\mathfrak{D}$ and $\mathfrak{D}$) for source injective
representation quivers.
\begin{cor}\label{quivgi}
Let $Q$ be a source injective representation quiver and $R$ be a
Gorenstein ring (see \cite{iwanaga}). A representation
${\mathcal{M}}$ of the quiver $Q$ is Gorenstein injective if and
only if for each vertex $v$ ${\mathcal{M}}(v)$ is a Gorenstein
injective $R$-module and the sequence
$$0\to \Ker(f_v)\to {\mathcal{M}}(v)\stackrel{f_v}{\to} \prod_{s(a)=v}{\mathcal{M}}(t(a))\to
0$$ is exact with $\Ker(f_v)$ Gorenstein injective (where
$f_v:{\mathcal{M}}(v)\to \prod_{s(a)=v}{\mathcal{M}}(t(a))$ is the
induced by the ${\mathcal{M}}(a):{\mathcal{M}}(v)\to
{\mathcal{M}}(t(a))$).
\end{cor}
{\bf Proof.} Suppose ${\mathcal{M}}$ is Gorenstein injective, then
there will exist an exact sequence of injectives
\begin{equation}\label{eq1}
\cdots \to {\mathcal{E}}_{-2}\to {\mathcal{E}}_{-1}\to
{\mathcal{E}}_0\to {\mathcal{E}}_1\to {\mathcal{E}}_2\to \cdots
\end{equation}
with ${\mathcal{M}}=\Ker({\mathcal{E}}_0\to {\mathcal{E}}_1)$, which
is $\Hom({\mathcal{E}},-)$ exact, for all injective representations.
Then, for a fixed vertex $v$ we have the corresponding exact
sequence of injective modules
$$\cdots \to {\mathcal{E}}_{-2}(v)\to {\mathcal{E}}_{-1}(v)\to {\mathcal{E}}_0(v)\to {\mathcal{E}}_1(v)\to
{\mathcal{E}}_2(v)\to \cdots$$ with
${\mathcal{M}}(v)=\Ker({\mathcal{E}}_0(v)\to {\mathcal{E}}_1(v))$.
Then, taking an integer sufficiently large in absolute value, and
apply \cite[Theorem 9.1.11(7)]{EdO}, we will have that
${\mathcal{M}}(v)$ is Gorenstein injective for all vertex $v$. Now,
since ${\mathcal{M}}$ is Gorenstein injective, there will exist a
short exact sequence of representations $$0\to {\mathcal{U}}\to
{\mathcal{E}}\to {\mathcal{M}}\to 0,$$ with ${\mathcal{E}}$
injective and ${\mathcal{U}}$ Gorenstein injective. By Proposition
\ref{epies}, ${\mathcal{E}}(v)\to \prod_{s(a)=v}{\mathcal{E}}(t(a))$
is a splitting epimorphism, so it easily follows that
${\mathcal{M}}(v)\to \prod_{s(a)=v}{\mathcal{M}}(t(a))$ is an
epimorphism. Finally let us show that $ \Ker(f_v) $ is Gorenstein
injective. By Proposition \ref{epies} we have, from (\ref{eq1}) the
exact sequence of injective kernels:
$$\cdots\to \Ker(t_v)^{-2}\to \Ker(t_v)^{-1}\to \Ker(t_v)^{0}\to \Ker(t_v)^{1}\to \Ker(t_v)^{2}\to \cdots  $$
where $$0\to \Ker(t_v)^{j}\to E_j(v)\stackrel{t_v}{\longrightarrow}
\prod_{s(a)=v}E_j(t(a))\to 0,$$ for all $j\in \Z$, such that
$\Ker(f_v)=\Ker(\Ker(t_v)^0\to \Ker(t_v)^1)$ which is $\Hom(E,-)$
exact for all injective $R$-modules $E$, by a reasoning similar to
the preceding.

Let us see that the conditions are sufficient. Since $(Q,R\Mod)$ is
a Gorenstein category we may find a short exact sequence
$$0\to {\mathcal{G}}\to {\mathcal{L}}\to {\mathcal{M}}\to 0,$$ with
${\mathcal{G}}$ a Gorenstein injective representation and
${\mathcal{L}}$ a representation of finite injective dimension. By
the previous proof ${\mathcal{G}}(v)$ is a Gorenstein injective
$R$-module and $\Ker(h_v)$ in the exact
$$0\to \Ker(h_v)\to {\mathcal{G}}(v)\stackrel{h_v}{\longrightarrow} \prod_{s(a)=v}{\mathcal{G}}(t(a))\to
0$$ is also Gorenstein injective. Moreover, since ${\mathcal{L}}$ is
of finite injective dimension, it is obvious that ${\mathcal{L}}(v)$
and the kernel of ${\mathcal{L}}(v)\to
\prod_{s(a)=v}{\mathcal{L}}(t(a))$ is also of finite injective
dimension (and this map is a surjection). So by the definition of a
source injective representation quiver it follows that
${\mathcal{L}}$ is an injective representation, so ${\mathcal{M}}$
will be also Gorenstein injective.{\hfill$\Box$}

\bigskip\par
In a dual manner and using the results of \cite{EnEst} we can
characterize Gorenstein projective representations of quivers.
\begin{cor}\label{quivgp}
Let $Q$ be a sink projective representation quiver and let $R$ be a
Gorenstein ring. A representation ${\mathcal{M}}$ of the quiver $Q$
is Gorenstein projective if and only if ${\mathcal{M}}(v)$ is a
Gorenstein projective $R$-module, the sequence
$$0\to \oplus_{t(a)=v} {\mathcal{M}}(s(a))\stackrel{f_v}{\longrightarrow} {\mathcal{M}}(v)\to {\mathcal{C}}(v)\to
0$$ is exact and ${\mathcal{C}}(v)$ is Gorenstein projective (where
$f_v:\oplus_{t(a)=v} {\mathcal{M}}(s(a))\to {\mathcal{M}}(v)$ is the
induced by ${\mathcal{M}}(s(a))\to {\mathcal{M}}(v)$).
\end{cor}

We can also characterize the class of all representations of finite
injective dimension for quivers in $\mathfrak{I}$.
\begin{prop}\label{nue}
A representation ${\mathcal{M}}$ of a source injective
representation quiver is of finite injective dimension if and only
if $\Sup\{injdim_R {\mathcal{M}}(v)<\infty:\ v\in V\}$ is finite.
\end{prop}
{\bf Proof.} It is obvious that the condition is necessary, for if
$injdim_Q {\mathcal{M}}\leq n$ then, for all vertex $v\in V$,
$injdim_R {\mathcal{M}}(v)\leq n$.

Conversely, suppose that $n=\Sup\{injdim_R {\mathcal{M}}(v)<\infty:\
v\in V\}$, and let $v\in V$ be a vertex, then there exists an exact
sequence of injective $R$-modules, $$0\to {\mathcal{M}}(v)\to
E_0^v\to E_1^v\to \cdots \to E_n^v\to 0$$ But then we have the short
exact sequence of representations
$$0\to {\mathcal{M}}\to \prod_{v\in V} e_*^v(E_0^v)\to {\mathcal{M}}_1\to 0$$
such that $injdim_R {\mathcal{M}}_1(v)\leq n-1$, for every $ v\in
V$. Now we repeat the same procedure to get ${\mathcal{M}}_2$ from
${\mathcal{M}}_1$ satisfying $injdim_R {\mathcal{M}}_2(v)\leq n-2$,
for all $v\in V$. So finally we get the exact
$$0\to {\mathcal{M}}\to {\mathcal{E}}_0\to {\mathcal{E}}_1\to \cdots \to {\mathcal{E}}_{n-1}\to{\mathcal{M}}_n\to 0$$
such that, ${\mathcal{M}}_n(v)$ is injective, for all vertex $ v\in
V$. Let us show that $injdim_Q M\leq 1$. For if $$0\to
{\mathcal{M}}_n\to {\mathcal{E}}_n\to {\mathcal{C}}\to 0$$ is an
exact sequence with ${\mathcal{E}}_n$ injective, by Proposition
\ref{epies}, we have that ${\mathcal{C}}(v)$ is an injective
$R$-module and ${\mathcal{C}}(v)\to
\prod_{s(a)=v}{\mathcal{C}}(t(a))$ is an epimorphism with injective
kernel, for all $v\in V$, so since the quiver is a source injective
representation quiver, ${\mathcal{C}}$ will be an injective
representation. {\hfill$\Box$}

\par\noindent
{\bf Remark.} Notice that the previous proposition shows that if
$\Sup\{injdim_R {\mathcal{M}}(v)<\infty:\ v\in V\}\leq n$ then
$injdim_Q{\mathcal{M}}\leq n+1$. The converse clearly is not true.

When $(Q,R\Mod)$ is a Gorenstein category, we can define the {\it
Gorenstein injective dimension} of a representation ${\mathcal{M}}$
as the least natural number $k$ such that the $n$th cosyzygy of
${\mathcal{M}}$ is Gorenstein injective, that is, $k$ is the least
integer such that there exists an exact sequence
$$0\to {\mathcal{M}}\to {\mathcal{E}}_0\to {\mathcal{E}}_1\to \cdots
\to {\mathcal{E}}_{k-1}\to {\mathcal{G}}_k\to 0$$ with
${\mathcal{E}}_i$ injective $0\leq i\leq k-1$ and ${\mathcal{G}}_k$
Gorenstein injective ($k=\infty$ if there exists no such a natural
number). We denote the Gorenstein injective dimension of a
representation by $Ginjdim_Q {\mathcal{M}}$.

\begin{prop}
Let $(Q,R\Mod)$ be a Gorenstein category and let ${\mathcal{M}}$ be
a representation. Suppose that ${\mathcal{M}}(v)$ is a Gorenstein
injective $R$-module for all $v\in V$, then $Ginjdim_Q
{\mathcal{M}}\leq 1$.
\end{prop}
{\bf Proof.} We consider a short exact sequence $0\to
{\mathcal{M}}\to {\mathcal{G}}\to {\mathcal{C}}\to 0$ with
${\mathcal{G}}$ a Gorenstein injective representation and $injdim_Q
{\mathcal{C}}<\infty$. Let ${\mathcal{E}}\to {\mathcal{G}}\to 0 $ be
exact with ${\mathcal{E}}$ injective. We have the pullback diagram
 \DIAGV{50} {}  \n{}  \n{0}\n{}\n{0}\nn
{}\n{}\n{\sar}\n{}\n{\sar}\nn
{}\n{}\n{{\mathcal{K}}}\n{\eeql}\n{{\mathcal{K}}}\nn
{}\n{}\n{\sar}\n{}\n{\sar}\nn{0}
\n{\ear}\n{{\mathcal{P}}}\n{\ear}\n{{\mathcal{E}}}\n{\ear}\n{{\mathcal{C}}}\n{\ear}\n{0}\nn
{}\n{}\n{\sar}\n{}\n{\sar}\n{}\n{\seql}\nn
{0}\n{\ear}\n{{\mathcal{M}}}\n{\ear}\n{{\mathcal{G}}}\n{\ear}\n{{\mathcal{C}}}\n{\ear}\n{0}\nn
{}\n{}\n{\sar}\n{}\n{\sar}\nn {}\n{}\n{0}\n{}\n{0}\diag Now
${\mathcal{K}}$ is a Gorenstein injective representation, so by
Corollary \ref{quivgi} and the hypothesis, we have that
${\mathcal{P}}(v)$ is a Gorenstein injective $R$-module. Since
$injdim_Q {\mathcal{C}}<\infty$ it follows that $injdim_Q
{\mathcal{P}}<\infty$, so by Proposition \ref{nue} we conclude that
${\mathcal{P}}(v)$ is an injective $R$-module. But then the previous
remark says that $injdim_Q {\mathcal{P}}\leq 1$. Now we form the
pushout diagram of $0\to {\mathcal{P}}\to
{\mathcal{E}}({\mathcal{P}})$ (the injective hull) and
${\mathcal{P}}\to {\mathcal{M}}\to 0$
 \DIAGV{55} {}  \n{}  \n{0}\n{}\n{0}\nn
{}\n{}\n{\sar}\n{}\n{\sar}\nn
{}\n{}\n{{\mathcal{K}}}\n{\eeql}\n{{\mathcal{K}}}\nn
{}\n{}\n{\sar}\n{}\n{\sar}\nn{0}
\n{\ear}\n{{\mathcal{P}}}\n{\ear}\n{{\mathcal{E}}({\mathcal{P}})}\n{\ear}\n{{\mathcal{T}}}\n{\ear}\n{0}\nn
{}\n{}\n{\sar}\n{}\n{\sar}\n{}\n{\seql}\nn
{0}\n{\ear}\n{{\mathcal{M}}}\n{\ear}\n{{\mathcal{Q}}}\n{\ear}\n{{\mathcal{T}}}\n{\ear}\n{0}\nn
{}\n{}\n{\sar}\n{}\n{\sar}\nn {}\n{}\n{0}\n{}\n{0}\diag Since
${\mathcal{K}}$ is Gorenstein injective, ${\mathcal{Q}}$ will be
Gorenstein injective and since $injdim_Q {\mathcal{P}}\leq 1$,
${\mathcal{T}}$ will be injective, so $Ginjdim_Q {\mathcal{M}}\leq
1$. {\hfill$\Box$}

\bigskip\par
{\bf Remark.} The argument of the previous proposition can be easily
extended and show that if $Ginjdim_Q {\mathcal{M}}(v)\leq k$, for
all $ v\in V$ then $Ginjdim_Q {\mathcal{M}}\leq k+1$.

In \cite{EnLoyBlas}, the authors characterize flat representations
of left rooted quivers. We may use this result with Theorem
\ref{injrooted} to relate flat and injective representations. We
need to introduce the following notation: we denote by
$(Q^{op},R\Mod)$ the category of representations over the quiver
$Q^{op}=(V,\Gamma^{op})$, with the same set of vertices as $Q$ and
where the arrows are the arrows of $\Gamma$ reversed (so $v\to w\in
\Gamma^{op}\Leftrightarrow w\to v\in \Gamma$). It is obvious that if
$Q$ is left rooted then $Q^{op}$ is right rooted.
\begin{cor}\label{flatinj}
Let $Q$ be a left rooted quiver. Then ${\mathcal{F}}$ is a flat
representation over $(Q,R\Mod)$ if, and only if, ${\mathcal{F}}^+$
is an injective representation over the quiver
$(Q^{op},Mod\mbox{-R})$ (where ${\mathcal{F}}^+$ is the
representation given by
${\mathcal{F}}^+(v)=\Hom_{\Z}({\mathcal{F}}(v),\Q/\Z)$ and, for an
arrow $a:v\to w$ and the morphism $f_a:{\mathcal{F}}(v)\to
{\mathcal{F}}(w)$,
$f^+_{a^{op}}=\Hom_{\Z}(f_a,\Q/\Z):{\mathcal{F}}^+(w)\to
{\mathcal{F}}^+(v)$)
\end{cor}
{\bf Proof.} This follows from \cite[Theorem 3.7]{EnLoyBlas} and
Theorem \ref{injrooted} by noticing that
$$0\to
\oplus_{t(a)=v}{\mathcal{F}}(s(a))\stackrel{f_v}{\longrightarrow}
{\mathcal{F}}(v)\stackrel{p_v}{\to}{\mathcal{C}}(v)\to 0 $$ is a
pure exact sequence of left $R$-modules if, and only if,
$$0\to {\mathcal{C}}^+(v)\stackrel{p_v^+}{\longrightarrow}{\mathcal{F}}^+(v)
\stackrel{f_v^+}{\longrightarrow}\prod_{s(a)=v}{\mathcal{F}}^+(t(a))\to
0$$ is a split short exact sequence of right
$R$-modules.{\hfill$\Box$}

\begin{defn}
An object $M$ of a monoidal Grothendieck category ${\mathcal A}$
(see, for example, \cite[pg.157]{maclane} for the definition of a
monoidal category) is said to be Gorenstein flat if there is an
exact sequence
$$\cdots \to F_{-3}\to F_{-2}\to F_{-1}\to F_0\to F_{1}\to F_{2}\to
\cdots$$ of flat objects such that $M=\Ker(F_0\to F_{1})$ and such
that the sequence is $E\otimes -$exact for every injective object
$E$ (i.e. the functor $E\otimes -$ leaves the sequence exact).
\end{defn}

\begin{lem}\label{implica1}
Let $Q$ be a left rooted quiver and $R$ be a Gorenstein ring. Then
if ${\mathcal{M}}$ is a Gorenstein flat representation of $Q$, the
representation ${\mathcal{M}}^+$ of the quiver $Q^{op}$ is
Gorenstein injective.
\end{lem}
{\bf Proof.} Since ${\mathcal{M}}$ is Gorenstein flat, there exists
an exact sequence
\begin{equation}\label{plares}
\cdots \to {\mathcal{F}}_{-3}\to {\mathcal{F}}_{-2}\to
{\mathcal{F}}_{-1}\to {\mathcal{F}}_0\to {\mathcal{F}}_{1}\to
{\mathcal{F}}_{2}\to \cdots
\end{equation} of flat
representations with ${\mathcal{M}}=\Ker({\mathcal{F}}_0\to
{\mathcal{F}}_1)$ which remains exact when applying
${\mathcal{E}}\otimes -$, for all injective representation in $(Q,
\mbox{ Mod-}R)$, but then by Corollary \ref{flatinj} we have the
exact sequence
$$\cdots \to {\mathcal{F}}_2^+\to {\mathcal{F}}_1^+\to {\mathcal{F}}_0^+\to {\mathcal{F}}_{-1}^+\to
{\mathcal{F}}_{-2}^+\to\cdots$$ of injective representations, such
that ${\mathcal{M}}=\Ker({\mathcal{F}}_{-1}^+\to
{\mathcal{F}}_{-2}^+)$. We see that the last exact sequence is
$\Hom({\mathcal{E}},-)$ exact, for all injective representation
${\mathcal{E}}$. So let ${\mathcal{E}}$ be an injective
representation. From (\ref{plares}) we have an exact
$$\cdots \to {\mathcal{E}}\otimes {\mathcal{F}}_{-2}\to {\mathcal{E}}\otimes
{\mathcal{F}}_{-1}\to {\mathcal{E}}\otimes {\mathcal{F}}_0\to
{\mathcal{E}}\otimes {\mathcal{F}}_{1}\to {\mathcal{E}}\otimes
{\mathcal{F}}_{2}\to \cdots$$so we have an exact sequence
$$\cdots \to ({\mathcal{E}}\otimes {\mathcal{F}}_{2})^+\to ({\mathcal{E}}\otimes {\mathcal{F}}_{1})^+\to
({\mathcal{E}}\otimes {\mathcal{F}}_{0})^+\to ({\mathcal{E}}\otimes
{\mathcal{F}}_{-1})^+\to ({\mathcal{E}}\otimes
{\mathcal{F}}_{-2})^+\to \cdots$$in $Q^{op}$. But the canonical
isomorphism
$$({\mathcal{E}}\otimes_{RQ} {\mathcal{F}}_i)^+\cong \Hom_{RQ}({\mathcal{E}},{\mathcal{F}}_i^+),\ i\in \Z$$ still
holds for rings without unit (as the case of $RQ$ whenever $Q$ has
an infinite number of vertices) so we have the exact sequence
$$\cdots\to \Hom({\mathcal{E}},{\mathcal{F}}_{2}^+)\to \Hom({\mathcal{E}},
{\mathcal{F}}_{1}^+)\to \Hom({\mathcal{E}},{\mathcal{F}}_0^+)\to
\Hom({\mathcal{E}},{\mathcal{F}}_{-1}^+)\to
\Hom({\mathcal{E}},{\mathcal{F}}_{-2}^+)\to \cdots$$ {\hfill$\Box$}

\begin{lem}\label{implica2}
Let $Q$ be a forest whose connected components are barren trees and
$_R R$ be noetherian. Then a representation ${\mathcal{E}}$ of $Q$
is injective if and only if ${\mathcal{E}}^+$ is flat (over
$Q^{op})$
\end{lem}
{\bf Proof.} We will use the characterization of flat
representations given in \cite[Theorem 3.7]{EnLoyBlas} and Corollary
\ref{inybarren}. Notice that, for all $v\in V$, ${\mathcal{E}}(v)$
is an injective $R$-module if and only if ${\mathcal{E}}^+(v)$ is a
flat $R$-module (see for example \cite[Corollary 3.2.17]{EdO}) and
furthermore the sequence
$$0\to {\mathcal{K}}(v)\to {\mathcal{E}}(v)\to \prod_{s(a)=v}{\mathcal{E}}(t(a))\to 0$$ is exact, with
${\mathcal{K}}(v)$ injective if and only if $$0\to
\oplus_{t(a)=v}{\mathcal{E}}^+(s(a))\to {\mathcal{E}}^+(v)\to
{\mathcal{K}}^+(v)\to 0$$ is an exact sequence with
${\mathcal{K}}^+(v)$ flat (so hence pure exact sequence) for all
$v\in V$. Notice that, since the connected components of $Q$ are
barren trees, then, for every vertex $v\in V$ the set $\{ a\in
\Gamma: \ s(a)=v\}$ is finite. {\hfill$\Box$}
\begin{prop}
Let $Q$ be a forest whose connected components are barren trees and
let $R$ be Gorenstein. A representation ${\mathcal{M}}$ of the
quiver $Q$ is Gorenstein flat if and only if for each vertex $v$,
${\mathcal{M}}(v)$ is a Gorenstein flat $R$-module and the sequence
$$0\to \oplus_{t(a)=v} {\mathcal{M}}(s(a))\stackrel{f_v}{\longrightarrow}
{\mathcal{M}}(v)\to {\mathcal{C}}(v)\to 0$$ is exact and
${\mathcal{C}}(v)$ is Gorenstein flat (where $f_v:\oplus_{t(a)=v}
{\mathcal{M}}(s(a))\to {\mathcal{M}}(v)$ is a morphism induced by
${\mathcal{M}}(s(a))\to {\mathcal{M}}(v)$).
\end{prop}
{\bf Proof.} Necessity. By Lemma \ref{implica1}, ${\mathcal{M}}^+$
is Gorenstein injective, so by Corollary \ref{quivgi}
${\mathcal{M}}^+(v)$ is a Gorenstein injective $R$-module for all $
v\in V$, but then ${\mathcal{M}}(v)$ is a Gorenstein flat
$R$-module. Furthermore since $$0\to {\mathcal{C}}^+(v)\to
{\mathcal{M}}^+(v)\to \prod_{s(a)=v}{\mathcal{M}}^+(t(a))\to 0$$ is
exact with ${\mathcal{C}}^+(v)$ Gorenstein injective, then $$0\to
\oplus_{t(a)=v}{\mathcal{M}}(s(a))\to {\mathcal{M}}(v)\to
{\mathcal{C}}(v)\to 0$$is exact with ${\mathcal{C}}(v)$ Gorenstein
flat, for every $ v\in V$.

Sufficiency. Since $R$ is Gorenstein the hypothesis is equivalent to
saying that ${\mathcal{M}}^+(v)$ is Gorenstein injective and that
$$0\to {\mathcal{C}}^+(v)\to {\mathcal{M}}^+(v)\to \prod_{s(a)=v}{\mathcal{M}}^+(t(a))\to 0$$ is exact
with ${\mathcal{C}}^+(v)$ Gorenstein injective. But, by Corollary
\ref{quivgi}, this implies that ${\mathcal{M}}^+$ is a Gorenstein
injective representation. So there exists an exact sequence
$${\mathcal{E}}_n\to {\mathcal{E}}_{n-1}\to \cdots {\mathcal{E}}_0\to {\mathcal{M}}^+\to 0$$
with ${\mathcal{E}}_i$ injective representations for all $i$ such
that $1\leq i\leq n$. Now, by Lemma \ref{implica2}, $$0\to
{\mathcal{M}}^{++}\to {\mathcal{E}}_0^+\to\cdots\to
{\mathcal{E}}_{n-1}^+\to {\mathcal{E}}_n^+$$ is an exact sequence of
flat representations. So $\Tor^i({\mathcal{L}},M^{++})=0$, for all
$i\geq 1$ and for all representation ${\mathcal{L}}$ such that
$injdim_Q {\mathcal{L}}<\infty$. But now there exists a pure exact
sequence
$$0\to {\mathcal{M}}\to {\mathcal{M}}^{++}\to {\mathcal{M}}^{++}/{\mathcal{M}}\to 0$$
of representations, so
$\Tor^1({\mathcal{L}},{\mathcal{M}}^{++}/{\mathcal{M}})=0$ for all
such ${\mathcal{L}}$. But then
$\Tor^i({\mathcal{L}},{\mathcal{M}}^{++}/{\mathcal{M}})=0$, for
every $i\geq 1$. So then $\Tor^i({\mathcal{L}},{\mathcal{M}})=0$ for
each $i\geq 1$ and for all representations ${\mathcal{L}}$ such that
$injdim_Q {\mathcal{L}}<\infty$. But this means that if we take a
complex
$$\cdots \to {\mathcal{F}}_{-2}\to {\mathcal{F}}_{-1}\to
{\mathcal{F}}_0\to {\mathcal{F}}_1\to {\mathcal{F}}_2\to \cdots$$ such
that ${\mathcal{M}}=\Ker({\mathcal{F}}_0\to {\mathcal{F}}_1)$ with
${\mathcal{F}}_i$ a flat representation for all $i\in\Z$, then this
will be ${\mathcal{E}}\otimes -$ exact for all injective
representation $\mathcal{E}$. Furthermore, by \cite[Theorem
3.7]{EnEsGar}, $injdim_{RQ} RQ<\infty$ so then the complex will be
exact. {\hfill$\Box$}

\bigskip\begin{center}
{\bf ACKNOWLEDGEMENTS}
\end{center}

The authors wish to thank the referee for many helpful comments and
suggestions.


\begin{thebibliography}{99}
\bibitem{assem}
I. Assem, D. Simson and A. Skowronski. {\sl Elements of the
representation theory of associative algebras.} Vol. 1. Techniques
of representation theory. London Mathematical Society Student Texts,
65. Cambridge University Press, Cambridge, 2006.

\bibitem{EnEst}
E. Enochs and S. Estrada. {\sl Projective representations of
quivers.} Comm. Algebra. {\bf 33}(2005), 3467-3478.
\bibitem{Estr}
E. Enochs and S. Estrada {\sl Relative homological algebra in the
category of quasi-coherent sheaves.} Adv. in Math. {\bf 194} (2005),
284-295.
\bibitem{gorenstein}
E. Enochs, S. Estrada and J.R. Garc\'{\i}a Rozas. {\sl Gorenstein
Categories. Tate Cohomology on Projective Schemes.} To appear in
Math. Nachr.
\bibitem{EnEsGar}
E. Enochs, S. Estrada, J.R. Garc\'{\i}a Rozas and A. Iacob. {\sl
Gorenstein representations of quivers.} To appear in Arch. Math.
(Basel).
\bibitem{EnLoyPark}
E. Enochs, J.R. Garc\'{\i}a Rozas, L. Oyonarte and S. Park. {\sl
Noetherian Quivers.} Quaest. Math. {\bf 25(4)} (2002), 531-538.
\bibitem{EnHer}
E. Enochs and I. Herzog. {\sl A homotopy of quiver morphisms with
applications to representations.} Canad. J. Math. {\bf 51(2)}
(1999), 294-308.
\bibitem{EdO}
E. Enochs and O. Jenda. {\sl Relative homological algebra.} De
Gruyter Expositions in Mathematics {\bf 30}, de Gruyter, Berlin,
2000.
\bibitem{EnLoyBlas}
E. Enochs, L. Oyonarte and B. Torrecillas. {\sl Flat covers and flat
representations of quivers.} Comm. Algebra. {\bf 32(4)}(2004),
1319-1338.
\bibitem{Estrad}
S. Estrada. {\sl Monomial algebras over infinite quivers.
Applications to $N$-complexes of modules.} To appear in Comm.
Algebra.
\bibitem{Gabriel}
P. Gabriel. {\sl Unzerlegbare Darstellungen I.} Manuscripta Math.
{\bf 6} (1972), 71-103.
\bibitem{hovey}
M. Hovey. {\sl Model Categories.} AMS, Mathematical Surveys and
Monographs, {\bf 63}, 1999.
\bibitem{iwanaga}
Y. Iwanaga. {\sl On rings with finite self-injective dimension.}
Comm. Algebra. {\bf 7(4)} (1979), 393-414.

\bibitem{procesi}
L. Le Bruyn and C. Procesi. {\sl Semisimple representations of
quivers.} Trans. Amer. Math. Soc. {\bf 317(2)} (1990), 585-598.

\bibitem{maclane}
S. Mac Lane. {\sl Categories for the working mathematician.}, 2nd.
ed., Graduate Texts in Mathematics {\bf 5}, Springer-Verlag, New
York, 1998.
\bibitem{sten}
B. Stenstr\"om. {\sl Rings of quotients.} Springer-Verlag, 1975.
\end{thebibliography}
\end{document}